\documentclass[12pt,twoside]{article}
\usepackage{amsmath,amssymb,mathrsfs,graphicx}
\usepackage[small,sc]{caption} 
\usepackage{xcolor}
\usepackage{tikz} 
\usetikzlibrary{decorations.pathreplacing}

\newtheorem{theorem}{Theorem}[section]
\newtheorem{lemma}[theorem]{Lemma}
\newtheorem{proposition}[theorem]{Proposition}
\newtheorem{corollary}[theorem]{Corollary}

\newtheorem{rmrk}[theorem]{Remark}

\DeclareMathAlphabet{\mathbfit}{OML}{cmm}{b}{it}

\makeatletter
\@addtoreset{equation}{section}
\makeatother

\newenvironment{remark}
{\begin{rmrk} \em}
{\end{rmrk}}

\newcommand{\fn} {function}
\newcommand{\me} {measure}

\newcommand{\erg} {ergodic}
\newcommand{\sy} {system}

\newcommand{\dsy} {dynamical system}

\newcommand{\R} {\mathbb{R}}

\newcommand{\Z} {\mathbb{Z}}
\newcommand{\N} {\mathbb{N}}
\newcommand{\qed} {\hfill {\small Q.E.D.} \par\medskip}
\newcommand{\skippar} {\par\medskip}
\newcommand{\ds} {\displaystyle}
\newcommand{\proof} {\noindent \textsc{Proof.} }
\newcommand{\proofof}[1] {\noindent \textsc{Proof of {#1}.} }
\newcommand{\article}[3] {\textsc{{#1}}, {\itshape {#2}}, {{#3}}.}
\newcommand{\book}[3] {\textsc{{#1}}, {\itshape {#2}}, {{#3}}.}
\newcommand{\vol} {\textbf}
\newcommand{\eps} {\varepsilon}

\newcommand{\rset}[2] {\left\{ #1 \: \left| \: #2 \right. \! \right\} }
\newcommand{\lset}[2] {\left\{ \left. \! #1 \: \right| \: #2 \right\} }

\newcommand{\into} {\longrightarrow}

\newcommand{\newfig}[4] {
\bigskip
\begin{figure}[htbp]
  \centering
  \includegraphics[width=#3]{#2}
  \begin{minipage}[t]{0.80\linewidth} 
    \caption{#4}
    \protect\label{#1}
  \end{minipage}
\end{figure}
}

\newcommand{\m} {mixing}
\newcommand{\ob} {observable}

\newcommand{\ps} {\mathcal{M}}   
\newcommand{\sca} {\mathscr{A}}   
\newcommand{\scb} {\mathscr{B}} 
\newcommand{\leb} {m}   
\newcommand{\scv} {\mathscr{V}}
\newcommand{\iv} {V \nearrow \ps}
\newcommand{\ivv} {\scv \ni V \nearrow \ps}
\newcommand{\ivvlim} {\lim_{\ivv}}
\newcommand{\go} {\mathcal{G}}   
\newcommand{\lo} {\mathcal{L}}   
\newcommand{\ui} {{(0,1]}}   
\newcommand{\rp} {{\R_0^+}}   
\newcommand{\br} {\tau}   
\newcommand{\ibr} {\phi}   
\newcommand{\avg} {\overline{\mu}}   
\newcommand{\avgv} {\overline{\mu}_\scv} 
\newcommand{\avgleb} {\overline{\leb}}   
\newcommand{\x} {\xi}   
\newcommand{\trop} {\widehat{T}}   
\newcommand{\symb} {\oplus}   
\newcommand{\gou} {\mathcal{G}_\mathrm{unif}}   

\begin{document}

\title{\textbf{Uniformly global observables for 1D maps with an indifferent 
  fixed point}}

\author{
\scshape
Giovanni Canestrari\thanks{
Scuola di Dottorato in Matematica, Universit\`a di Roma Tor Vergata, Via 
della Ricerca Scientifica 1, 00133 Roma, Italy. E-mail: 
\texttt{canestrari@mat.uniroma2.it}.}
\ and
Marco Lenci\thanks{
Dipartimento di Fisica e Astronomia, Universit\`a di Bologna, 
Via Irnerio 46, 40126 Bologna, Italy. 
E-mail: {\texttt{marco.lenci@unibo.it}}.}
\thanks{
Istituto Nazionale di Fisica Nucleare,
Sezione di Bologna, Viale Berti Pichat 6/2,
40127 Bologna, Italy.}
}

\date{May 2024}

\maketitle

\begin{abstract}
  We study the property of global-local mixing for full-branched expanding maps 
  of either the half-line or the interval, with one indifferent fixed point. 
  Global-local mixing expresses the decorrelation of global vs local observables 
  w.r.t.\ to an infinite measure $\mu$. Global observables are essentially bounded 
  functions admitting an infinite-volume average, i.e., a limit for the average of the 
  function over bigger and bigger intervals; local observables are integrable 
  functions (both notions are relative to $\mu$). Of course, the definition of global
  observable depends on the exact definition of infinite-volume average. The first 
  choice for it would be to consider averages over the entire space minus a 
  neighborhood of the indifferent fixed point (a.k.a.\ the ``point at infinity''), in the 
  limit where such neighborhood vanishes. This is the choice that was made in 
  previous papers on the subject. The classes of systems for which global-local
  mixing was proved, with this natural choice of global observables, are ample 
  but not really general. In this paper we consider \emph{uniformly global 
  observables}, i.e., $L^\infty$ functions whose averages over \emph{any}
  interval $V$ converges to a limit, uniformly as $\mu(V) \to \infty$. Uniformly 
  global observables form quite an extensive subclass of all global observables. 
  We prove global-local mixing in the sense of uniformly global observables, for 
  two truly general classes of expanding maps with one indifferent fixed point, 
  respectively on $\mathbb{R}_0^+$ and on $(0,1]$. The technical core of the 
  proofs is rather different from previous work.

  \bigskip\noindent 
  \textbf{Mathematics Subject Classification (2020):} 37A40, 37A25, 37E05, 
  37D25, 37C25. 
  
  \bigskip\noindent
  \textbf{Keywords:} infinite ergodic theory, global-local mixing, infinite-volume 
  mixing, one-dimensional maps, intermittent maps, non-uniformly expanding 
  maps.
\end{abstract}

\section{Introduction}
\label{sec-intro}

In the context of infinite ergodic theory, the notion of \emph{global-local \m} 
involves two types of \ob\ functions: local \ob s and global \ob s. 

If $(\ps, \sca, \mu, T)$ is a \dsy\ where $(\ps, \sca, \mu)$ is a $\sigma$-finite, 
infinite measure space and $T: \ps \into \ps$ is surjective, bi-measurable, and 
non-singular w.r.t.\ $\mu$, a local \ob\ is a function $f \in L^1(\ps, \sca, \mu)$. 

The notion of global \ob\ requires a definition of infinite-volume average, which
is given as follows. Let $\scv$ be an \emph{exhaustive family}, that is, a 
collection of sets $V \in \sca$, with $\mu(V) < \infty$, containing at least an 
increasing sequence $(V_k)_{k \in \N}$ such that $\bigcup_k V_k = \ps$. The
\emph{infinite-volume limit} $\ivv$ is defined to be the uniform limit for $V \in \scv$,
as $\mu(V) \to \infty$. In particular, a function $F: \ps \into \R$ is said to have
\emph{infinite-volume average} $\avgv(F)$, w.r.t.\ $\scv$ and $\mu$, if the 
following limit
\begin{equation} \label{def-avg}
  \avgv(F) := \ivvlim \, \frac1 {\mu(V)} \int_V F \, d\mu
\end{equation}
exists, which means that
\begin{equation} \label{def-avg2}
  \lim_{r \to \infty} \, \sup_{V \in \scv \atop \mu(V) \ge r} \, \left| \frac1 {\mu(V)} 
  \int_V F \, d\mu - \avgv(F) \right| = 0.
\end{equation}
A global \ob\ is an essentially bounded function admitting an infinite-volume
average, that is, an element of
\begin{equation} \label{gmax}
  \go_\scv(\ps, \sca, \mu) := \lset{F \in L^\infty(\ps, \sca, \mu)} {\exists \avgv(F)} .
\end{equation}
This space will be henceforth referred to as the \emph{(maximal) space of 
global \ob s}, relative to $\scv$ and $\mu$.%
\footnote{For a ``physical'' interpretation of global and local \ob s see 
\cite{lcmp, liutam, lpmu, dn}.}

\paragraph*{Convention.} From now on we will abbreviate all notation such as
$\go_\scv(\ps, \sca, \mu)$, $L^\infty(\ps, \sca, \mu)$, etc.\ to $\go_\scv$, 
$L^\infty$, etc. Should certain dependences need to be specified we might 
write $\go_\scv(\mu)$, $L^\infty(\mu)$, etc.
\bigskip\smallskip

The definitions of infinite-volume \m\ \cite{lcmp, lpmu} concern the decorrelation 
properties of global and local \ob s. Given subspaces $\go \subseteq 
\go_\scv$, $\lo \subseteq L^1$, we say that $T$ (endowed with the \me\ $\mu$) 
is \emph{global-local \m} w.r.t.\ $\scv$, $\go$ and $\lo$ if 
\begin{equation} \label{glm}
  \forall F \in \go, \, \forall g \in \lo, \quad \lim_{n \to \infty} \, 
  \mu((F \circ T^n) g) = \avgv(F) \, \mu(g),
\end{equation}
with the standard notation $\mu(g) := \int_\ps g\, d\mu$. The reasons 
one may want/need to restrict the class of global and local \ob s, respectively, 
to strict subspaces of $\go_\scv$ and $L^1$ has to do with the truth and 
feasibility of the sought result. In general, if possible, one would like to show 
that $T$ is global-local \m\ w.r.t.\ $\scv$, $\go_\scv$ and $L^1$. We refer to 
this property as \emph{full global-local \m}, relative to $\scv$.

An equivalent way to introduce definition (\ref{glm}) is as follows:
\begin{equation}
  \forall F \in \go, \, \forall g \in \lo, \, g \ge 0, \, \mu(g)=1, \quad \lim_{n \to \infty} 
  \, T_*^n \mu_g(F) = \avgv(F),
\end{equation}
where $\mu_g$ is the measure determined by $\frac{d\mu_g}{d\mu} = g$ 
and $T_*$ is the 
push-forward operator for $T$ acting on measures. This formulation makes 
it apparent that global-local \m\ defines a sort of ``convergence to equilibrium'' 
for a class of probability measures, where the infinite-volume average plays 
the role of the limit measure and the global \ob s play the role of 
test functions. This seems reasonable in all cases where the dynamics
tends to spread trajectories all over an infinite-measure space, and so the
evolution of probability measures cannot have a limit in any standard 
measure-theoretic way \cite{lcmp, lpmu, lunif, bgl1, bgl2, bl, dn, ghr}.

Global \ob s are also used in the context of \emph{global-global \m}, another 
notion of `infinite-volume \m' that was devised in \cite{lcmp} (see also
 \cite{liutam, lunif, dn, tz1, tz2}), and for the question of pointwise convergence
of Birkhoff averages in infinite \me\ \cite{lm, dln, bs}. In infinite \erg\ theory,
there exist other definitions related to \m\  \cite{s, ds}, including the so-called 
\emph{Hopf-Krickeberg \m} (a short list of recent papers in this area
include \cite{mt, p, pt, ghr2}, see also references therein).

Expanding maps with a sufficiently flat indifferent fixed point are paradigmatic 
examples of \dsy s with an infinite \me. For the reasons set forth earlier, they 
are generally expected to be global-local \m\ (though not global-global \m, as 
their ``\m\ region'' is essentially of finite measure in an infinite-measure space 
\cite{bgl1}). 

In this paper we consider full-branched expanding maps defined either on
$\ui$ or on $\rp$, with a single indifferent fixed point respectively in 0
or at $\infty$ (the latter case meaning that $T(x) \sim x$, as $x \to +\infty$). 
Their infinite-volume \m\ properties were studied in \cite{bgl1, bl} (and partly
in \cite{bgl2}). In \cite{bgl1} global-local \m\ was proved for 
a class of maps of the interval with an increasing and a decreasing branch 
(plus certain generalizations) and for a class of maps of the half-line with 
finitely or infinitely many increasing branches, which preserve the Lebesgue 
measure (plus certain generalizations). In \cite{bl} global-local \m\ was proved 
for two classes of maps, of the interval and the half-line respectively, with 
finitely or infinitely many increasing branches: in the case of the interval, the 
most relevant assumption was a certain growth condition on the branches of 
the map; in the case of the half-line, the analogue of such condition was still 
assumed, though it was much more natural, and no assumption was made on 
the preservation of the Lebesgue measure. The maps studied in \cite{bgl1} 
include the Fairy map and the ones studied in \cite{bl} include the classical 
Pomeau-Manneville maps and the Liverani-Saussol-Vaienti maps.

For all the above maps, the authors proved full global-local \m\ relative to 
the exhaustive families
\begin{equation}
  \scv_\mathrm{min} := 
  \begin{cases}
      \lset{[a,1]}{0 < a < 1}, & \text{for the case } \ps = \ui; \\
      \lset{[0,a]}{a > 0}, & \text{for the case } \ps = \rp.
    \end{cases}
\end{equation}
In both cases, $\scv_\mathrm{min}$ is essentially the smallest exhaustive 
family one can take, if it is to include an increasing sequence of 
finite-measure sets covering $\ps$ (whence the notation). Therefore 
$\go_{\scv_\mathrm{min}}$ is essentially the largest class of global 
\ob s within the framework presented earlier. Since $L^1$ is also the 
largest class of local \ob s within the theory, full global-local \m\ w.r.t.\
$\scv_\mathrm{min}$ is the most general result of its type. 

That said, the state of the art is not quite satisfactory, as the classes of maps
investigated so far are not as general as one expects global-local \m\ to hold 
for. Here we prove global-local \m\ for two general classes of maps, on $\ui$ 
and $\rp$, respectively, admitting increasing and decreasing full branches, 
though we slightly restrict the class of global \ob s. We consider in fact 
\emph{uniformly global} \ob s, that is, global \ob s relative to the exhaustive 
families 
\begin{equation} \label{v-unif}
  \scv_\mathrm{unif} :=
  \begin{cases}
      \lset{[a,b]}{0 < a < b \le 1}, & \text{for the case } \ps = \ui; \\
      \lset{[a,b]}{0 \le a < b}, & \text{for the case } \ps = \rp.
    \end{cases}
\end{equation}
In view of definitions (\ref{def-avg})-(\ref{gmax}), this means that averaging one 
of these \ob s over \emph{any} interval with large enough measure will result 
in a uniformly good approximation of the infinite-volume average. 
Understandably, $\go_{\scv_\mathrm{unif}}$ is a strict subspace of 
$\go_{\scv_\mathrm{min}}$. For example, in the case $\ps = \rp$ and $\mu = 
\leb$, the Lebesgue \me, it is easy to check that
\begin{equation}
  \lim_{r \to +\infty} \, \frac1r \int_0^r \cos \sqrt{x} \, dx = 0,
\end{equation}
while, for all $r>0$,
\begin{equation}
 \underset{a \in \rp}{\sup / \inf} \ \frac1r \int_a^{a+r} \cos \sqrt{x} \, dx = \pm 1,
\end{equation}
implying that $F(x) := \cos \sqrt{x}$ belongs to $\go_{\scv_\mathrm{min}}$ but 
not to $\go_{\scv_\mathrm{unif}}$.

Uniformly global \ob s have been previously studied, together with other 
classes of global \ob s, in more general contexts than the present one, cf.\
\cite{dn, dln}.

For the classes of maps we consider, cf.\ Sections \ref{subs-maps-rp}
and \ref{subs-maps-ui}, we prove full global-local \m\ for uniformly global 
\ob s. The 
technique of the proofs is different from what was used so far: the main 
tools in \cite{bgl1, bgl2, bl} were the so-called \emph{persistently 
monotonic} local \ob s, that is, monotonic $L^1$ functions which maintain 
their monotonicity when acted upon by the dynamics (namely, by the 
corresponding transfer operator). Persistently monotonic local \ob s form an 
invariant convex cone in $L^1$. Here we also use dynamics-invariant cones of 
local observables, but of a different nature, related to the logarithmic derivative 
of the \ob s.

\skippar

Here is an outline of the paper: In Section \ref{sec-setup} we first introduce the
class of uniformly global \ob s, for both the case of the half-line and of the unit 
interval, discussing the dependence of their infinite-volume average on the 
\me\ defined on the system. Then we define our maps and state the main
results, most notably the full global-local \m\ for uniformly global \ob s, w.r.t.\ 
relevant equivalence classes of \me s. Section \ref{sec-proofs} contains the 
core proofs.

\paragraph*{Acknowledgments.} This research was partially supported by 
the PRIN Grant \linebreak
2017S35EHN of the Ministry of University and Research (MUR), 
Italy. It is also part of the authors' activity within the DinAmicI community 
(\texttt{www.dinamici.org}) and within the Gruppo Nazionale di Fisica
Matematica, INdAM, Italy. C.~G.~was also partially supported by the MUR 
Excellence Department Projects MatMod@TOV and Math@TOV, awarded to the
Department of Mathematics, Universit\`a di Roma Tor Vergata. Both authors 
acknowledge the hospitality of the Simons Center for Geometry and Physics at 
Stony Brook University, where part of this research was carried out.

\section{Setup}
\label{sec-setup}

For the rest of the paper we will lighten the notation and write $\scv := 
\scv_\mathrm{unif}$ and $\gou(\mu) := \go_{\scv_\mathrm{unif}}(\ps, \scb, \mu)$, 
cf.\ (\ref{v-unif}) and (\ref{gmax}), where $(\ps, \scb)$ is either $\rp$ or $\ui$, 
endowed with its respective Borel $\sigma$-algebra, and $\mu$ is an
infinite, locally finite measure there. In other words, $\gou(\mu)$ is in either 
case the class of \emph{all} uniformly global \ob s, that is, essentially 
bounded \fn s $F: \ps \into \R$ for which there exists
$\avg(F) := \avg_{\scv_\mathrm{unif}} (F)$ such that
\begin{equation} \label{def-avg3}
  \lim_{r \to \infty} \, \sup_{[a,b] \in \scv_\mu(r)} \, \left| \frac1r
  \int_a^b F \, d\mu - \avg(F) \right| = 0,
\end{equation}
having introduced the notation 
\begin{equation} \label{v-mu-r}
  \scv_\mu(r) := \rset{[a,b] \in \scv} {\mu([a,b]) = r}.
\end{equation}
(Notice that (\ref{def-avg3})-(\ref{v-mu-r}) are equivalent to (\ref{def-avg2}).)
The Lebesgue \me, in whichever space, will always be denoted by $\leb$.

\subsection{Dependence of the infinite-volume average on the global 
observable and on the measure}

Let us make some simple observations on the dependence of $\avg(F)$
on $F$ and $\mu$. Setting
\begin{equation} 
  \omega :=
  \begin{cases}
      0, & \text{for the case } \ps = \ui; \\
      +\infty, & \text{for the case } \ps = \rp,
    \end{cases}
\end{equation}
we see that $\omega$ is the only ``point at infinity'' for the \me\ $\mu$, in 
the sense that the complement of every neighborhood of $\omega$ has
finite \me\ (by the assumptions on $\mu$). It is not hard to see that only the 
values of $F$ around $\omega$ count for its infinite-volume average:

\begin{proposition} \label{prop1}
  If $F \in \gou(\mu)$, $G \in L^\infty(\ps, \mu)$ and $G(x) - F(x) \to 0$, as
  $x \to \omega$, then $G \in \gou(\mu)$ and $\avg(G) = \avg(F)$.
\end{proposition}

To prove this proposition it will be convenient to have a notation for the
neighborhoods of $\omega$, which will be used in other parts of the 
paper as well. For $M \in \ps$, let
\begin{equation} 
  U_M :=
  \begin{cases}
      (0,M), & \text{for the case } \ps = \ui; \\
      (M,+\infty), & \text{for the case } \ps = \rp.
    \end{cases}
\end{equation}
\smallskip

\proofof{Proposition \ref{prop1}} By possibly using $F- \avg(F)1$ in lieu of
$F$, where $1(x) \equiv 1$, we can always assume $\avg(F) = 0$. Given 
$\eps>0$, fix $M$ such that $|G-F| \le \eps/3$ in $U_M$. Then, for all $r>0$
and $V \in \scv_\mu(r)$,
\begin{equation} 
\begin{split}
  \left| \int_V G \, d\mu \right| &\le \int_{V \setminus U_M} |G| \, d\mu + 
  \left| \int_{V \cap U_M} G \, d\mu \right| \\
  &\le \| G \|_\infty \, \mu(V \setminus U_M) + \left| \int_{V \cap U_M} F \, d\mu 
  \right| + \frac \eps 3 \, \mu(V \cap U_M) \\
  &\le \left( \| G \|_\infty + \| F \|_\infty \right) \mu(V \setminus U_M) + 
  \left| \int_V F \, d\mu \right| + \frac \eps 3 \, r.  
\end{split}
\end{equation}
Now, dividing by $r$ and taking the sup for $V \in \scv_\mu(r)$, we see that,
for all sufficiently large $r$, the first term of the above r.h.s.\ can be made 
smaller or equal to $\eps/3$ (since $ \mu(V \setminus U_M) \le 
\mu(\ps \setminus U_M) < \infty$), and the same for the second term 
(since $\avg(F) = 0$, see (\ref{def-avg3})).
\qed

Focusing instead on the dependence of the infinite-volume average on the
\me, we have:

\begin{proposition} \label{prop2}
  If $\mu, \nu$ are two infinite, locally finite \me s on $\ps$ with $\nu \ll \mu$
  and $\frac{d\nu}{d\mu}(x) \to c \ne 0$, for $x \to \omega$, then 
  $\overline{\nu} = \avg$, in the sense that $\gou(\nu) = \gou(\mu)$ 
  and $\overline{\nu}(F) = \avg(F)$ for all $F \in \gou(\nu)$.
\end{proposition}

\proof Let us first show that $\avg$ is a \emph{restriction} of 
$\overline{\nu}$, meaning that $\gou(\mu) \subseteq \gou(\nu)$ and, for all 
$F \in \gou(\mu)$, $\overline{\nu}(F) = \avg(F)$. Fix $F \in \gou(\mu)$. By 
Proposition \ref{prop1}, as $\iv$,
\begin{equation} 
  \int_V F \, d\nu = \mu(V) \, \frac1 {\mu(V)} \int_V F \, \frac{d\nu}{d\mu} \, d\mu 
  \sim \mu(V) c \, \avg(F), \\
\end{equation}
where $\sim$ denotes \emph{exact} asymptotics. Using $1$ in place of $F$ in 
the above, one also has that $\nu(V) \sim \mu(V) c$. Since the limit of the ratio 
equals the ratio of the limits for our uniform limit $\iv$ as well, our initial claim
is proved.

As for the proof that $\overline{\nu}$ is a restriction of $\avg$, one uses
that $\frac{d\nu}{d\mu} > 0$ in a neighborhood $U_M$ of $\omega$, where
\begin{equation} 
  \frac{d\mu}{d\nu}(x) = \left( \frac{d\nu}{d\mu}(x) \right)^{-1} \to c^{-1}, 
\end{equation} 
as $x \to \omega$. So one applies the first part of this proof with $U_M, 
\nu |_{U_M}, \mu |_{U_M}$ in place of $\ps, \mu, \nu$, respectively.
The sought assertion follows from the fact that, by Proposition \ref{prop1},
the infinite-volume averages of $F 1_{U_M}$ and $F$, relative to either 
$\mu$ or $\nu$, exist together and are the same (here $1_{U_M}$ is the 
indicator \fn\ of $U_M$).
\qed

\begin{remark}
  The above assertion was already proved for general global \ob s, i.e., 
  \fn s in $\go_{\scv_\mathrm{min}}$, cf.\ \cite[Rmk.~2]{bl}, but that result does
  not imply Proposition \ref{prop2}.
\end{remark}

Proposition \ref{prop2} is actually implied by a more complicated result which
gives fairly general sufficient conditions for $\avg$ to be a restriction of 
$\overline{\nu}$:

\begin{proposition} \label{prop:restriction}
  Let $\mu, \nu$ be two infinite, locally finite \me s on $\ps$ with $\mu \approx 
  \nu$ (i.e., $\mu$ and $\nu$ are mutually absolutely continuous). If, for 
  all sufficiently large $r > 0$,
  \begin{equation} \label{eq:flat}
    \lim_{M \to \omega} \, \sup_{\substack{V \in \scv_\mu(r) \\ V \subseteq U_{M}} } 
    \left[\sup_V \frac{d\nu}{d\mu} \left( \inf_V \frac{d\nu}{d\mu}\right)^{-1} - 1 \right]
    = 0,
  \end{equation}
  and
  \begin{equation} \label{eq:comparability}
    \sup_{V \in \scv_\mu(r)} \nu(V) < \infty,
  \end{equation}
  then $\avg$ is a restriction of $\overline{\nu}$.
\end{proposition}

The proof of this proposition is given in Section \ref{sec-proofs}.

\begin{remark} \label{rm:simpl-cond}
In the case where $\mu \approx \leb$, $\nu \approx \leb$, with $\frac{d\mu}{dm}, 
\frac{d\nu}{dm} \in C^1$, the assumption (\ref{eq:flat}) of Proposition 
\ref{prop:restriction} is implied by the following condition, that is somehow 
easier to check in applications:
\begin{equation} \label{eq:flat-baby}
  \lim_{x \to \omega} \frac{ \frac d {dx}\! \left( \frac{d\nu}{d\mu} \right) \! (x)} 
  {\left( \frac{d\nu}{d\leb} \right) \! (x)} = 0.
\end{equation}

In fact, for any given $r>0$ and $M \in \ps$, let us
consider $V \in \scv_{\mu}(r)$, $V \subset U_M$, as in (\ref{eq:flat}), and let 
$x_\mathrm{max}$ and $x_\mathrm{min}$ denote two points in $V$ where 
$\frac{d\nu}{d\mu}$ attains its maximum and minimum, respectively. Denoting 
the ordinary derivative with the usual apostrophe, we have 
\begin{equation} \label{eq:log-eq}
\begin{split}
  \log \sup_V \frac{d\nu}{d\mu} - \log \inf_V \frac{d\nu}{d\mu} &= 
  \int_{x_\mathrm{min}}^{x_\mathrm{max}} \left( \frac{d\nu}{d\mu} \right)' \left(
  \frac{d\nu}{d\mu} \right)^{-1} d\leb \\
  &\le \int_V \left| \left( \frac{d\nu}{d\mu} \right)' \right| \left(\frac{d\nu}{d\mu} 
  \frac{d\mu}{dm}\right)^{-1} d\mu \\
  &\le \sup_{U_M} \left[ \left| \left( \frac{d\nu}{d\mu}\right)' \right| \left( 
  \frac{d\nu}{d\leb}\right)^{-1}\right] r,
\end{split} 
\end{equation}
which, for fixed $r$, vanishes as $M \to \omega$, by (\ref{eq:flat-baby}).
The uniformity in $V$ of the above limit readily implies (\ref{eq:flat}).
\end{remark}

The previous propositions apply to the question of global-local \m\ by 
means of the following simple lemma, which applies to a general \dsy,
as introduced in Section \ref{sec-intro}.

\begin{lemma} \label{lem-primo-ultimo}
  Suppose that $\nu \ll \mu$ are two \me s on $\ps$ and $\go$ is a space of 
  global \ob s, for the exhaustive family $\scv$, where $\overline{\nu}_\scv$ 
  and $\overline{\mu}_\scv$ coincide. If $(T, \mu)$ is global-local \m\ w.r.t.\ 
  $\scv$, $\go$ and $L^1(\mu)$, then $(T,\nu)$ is global-local \m\ w.r.t.\ $\scv$, 
  $\go$ and $L^1(\nu)$.
\end{lemma}

\proof Take $F \in \go$ and $g \in L^1(\nu)$. Clearly, $g \frac{d\nu}{d\mu}
\in L^1(\mu)$. So, for $n \to \infty$,
\begin{equation}
  \nu((F \circ T^n)g) = \mu \!\left( (F \circ T^n) \, g \frac{d\nu}{d\mu} \right) 
  \ \to \ \overline{\mu}_\scv(F) \, \mu \!\left( g \frac{d\nu}{d\mu} \right) =
  \overline{\nu}_\scv(F) \, \nu(g).
\end{equation}
\qed

\skippar

We are going to state our main theorem separately for the cases $\ps = \rp$ 
and $\ps = \ui$. It will be morally the same result in both cases, but the two 
theorems will not be equivalent. Obviously, by means of a suitable conjugation 
$\Psi: \ui \into \rp$ (say a smooth, decreasing \fn\ with a non-integrable singularity 
in 0), one can always view a \dsy\ $(\rp, \scb_{\rp}, \mu, T)$ as the \sy\ 
$(\ui, \scb_{\ui}, \mu_o, T_o)$, with $\mu_o := \Psi_*^{-1} \mu$ and $T_o := 
\psi^{-1} \circ T \circ \psi$, whence a theorem for the former \sy\ can be rewritten 
as a theorem for the latter \sy, and viceversa. But the assumptions on the map in
one case can become quite cumbersome when ``translated'' in terms of its 
conjugate (this occurs, for instance, for the expansivity of the map, which is not 
necessarily preserved by a conjugation of the above type); hence the 
convenience of stating a different theorem for each case.

\subsection{Maps of the half-line}
\label{subs-maps-rp}

Let $T : \rp \into \rp$ a Markov map w.r.t.\ the partition $\{ I_j \}_{j=0}^{N-1}$,
where $I_0 = (a_1, +\infty)$ and $I_j = (a_{j+1}, a_j)$, for $j=1, \ldots,, N-1$.
Here $0 = a_N < a_{N-1} < \cdots < a_1 \in \rp$. Clearly, $\{ I_j \}$ is a partition 
only up to a Lebesgue-null set, a subtlety that we shall henceforth 
forget. Assume the following:

\begin{itemize} 
\item[(A1)] $T |_{I_j}$ is a bijective map onto $\R^+$, and has an extension 
$\br_j$ which is defined on $[a_1, +\infty)$ (for $j=0$), or $[a_{j+1}, a_j)$ or 
$(a_{j+1}, a_j]$ (for $j \ge 1$, depending on $\br_j$ being increasing or 
decreasing, respectively). The \emph{branch} $\br_j$ is $C^2$ (this means
up to the closed endpoint of its domain).

\item[(A2)] There exists $\Lambda > 1$ such that $| \br'_j | \ge \Lambda$, 
for all $j \ge 1$. Also, $\br'_0 > 1$ and $\ds \lim_{x \to +\infty} \br_0'(x) = 1$.

\item[(A3)] $\ds \lim_{x \to +\infty} \,
\frac{ \br_0''(x) } { \br_0'(x)-1 } = 0$.

\item[(A4)] For all $j \ge 1$, either $\ds \lim_{x \to a_j^-} 
\frac{ \br_j''(x) } {( \br_j'(x) )^2} = 0$ or $\ds \lim_{x \to a_{j+1}^+} 
\frac{ \br_j''(x) } {( \br_j'(x) )^2} = 0$, depending on $\br_j$ being increasing or
decreasing, respectively.

\item[(A5)] $T$ is exact w.r.t.\ $\leb$, the Lebesgue \me\
on $\rp$.
\end{itemize}

\begin{remark}
  In \cite[App.~A]{bgl1} it is shown that, under the assumptions (A1)-(A4) 
  and
  \begin{itemize}
  \item[(A5)'] The \fn\ $u(x) := x - \br_0 (x)$ is positive, convex and 
    vanishing (hence decreasing), as $x \to +\infty$. Furthermore, 
    $u''$ is decreasing (hence vanishing),
  \end{itemize}
  $T$ is exact w.r.t.\ $\leb$. (In fact, Theorem A.1 of \cite{bgl1} proves that,
  under such assumptions and the existence of an invariant \me\ $\mu$,
  mutually absolutely continuous with $\leb$, $T$ is conservative and exact.
  But this last assumption is only needed for the conservativity of $T$.)
\end{remark}

\begin{theorem} \label{thm-main-rp}
  A map $T$ satisfying (A1)-(A5) is fully global-local mixing, relative to $\leb$, 
  w.r.t.\ $\scv = \scv_\mathrm{unif}$.
\end{theorem}

For the convenience of the reader we recall that the above statement
means that $T$ is global-local mixing for the Lebesgue \me, relative to all 
uniformly global \ob s and all integrable local \ob s. An immediate 
consequence of Theorem \ref{thm-main-rp} is that $T$ is global-local mixing for 
a large class of other \me s, though not fully, but relative to a large subclass of 
uniformly global \ob s (and all integrable local \ob s).

\begin{corollary} \label{cor-main-rp}
  Let $q \in (0,1]$ and $\nu \ll \leb$ on $\rp$ with $\frac{d\nu}{d\leb}$ bounded
  and $\frac{d\nu}{d\leb}(x) \sim x^{-q}$, as $x \to +\infty$. Under the same 
  hypotheses as Theorem \ref{thm-main-rp}, $\avgleb$ is a strict restriction of 
  $\overline{\nu}$ and $(\rp, \scb, \nu, T)$ is global-local mixing w.r.t.\ $\scv$, 
  $\gou(\leb)$ and $L^1(\nu)$.
\end{corollary}

\proofof{Corollary \ref{cor-main-rp}} Let us denote by $\lambda_q$ the \me\ with 
$\frac{d\lambda_q}{d\leb}(x) = (x+1)^{-q}$. By Proposition \ref{prop:restriction} 
and Remark \ref{rm:simpl-cond}, $\avgleb$ is a restriction of 
$\overline{\lambda_q}$. On the other hand, by Proposition \ref{prop2}, 
$\overline{\nu} = \overline{\lambda_q}$, so $\avgleb$ is a restriction of 
$\overline{\nu}$. To show that it is a strict restriction, it suffices to produce a 
global \ob\ in $\gou(\lambda_q) \setminus \gou(\leb)$. For $q \in (0,1)$, the \fn\ 
$F_q(x) := \cos( (x+1)^{1-q} )$ does the job because, for all $a \in \rp$ and 
$r > 0$,
\begin{equation}
  \frac1r \int_a^{a+r} F_q \, d\lambda_q = \frac{\sin( (a+r+1)^{1-q} ) - 
  \sin( (a+1)^{1-q} )} {r(1-q)},
\end{equation}
implying that $F_q \in \gou(\lambda_q)$ with $\overline{\lambda_q}(F_q) = 0$. 
For $q = 1$, one can take instead $F_1(x) := \cos(\log(x+1))$. For all $a \in \rp$
and $r > 0$,
\begin{equation}
  \frac1r \int_a^{a+r} F_1 \, d\lambda_1 = \frac{\sin( \log(a+r+1) ) - 
  \sin( \log(a+1))} {r},
\end{equation}
implying again that $F_1 \in \gou(\lambda_1)$ with 
$\overline{\lambda_1}(F_1) = 0$. On the other hand, for every $q \in (0,1]$ and 
$r > 0$, clearly
\begin{equation}
  \underset{a \to +\infty}{\limsup / \liminf} \ \frac1r \int_a^{a+r} F_q \, d\leb = \pm 1,
\end{equation}
whence $F_q \not\in \gou(\leb)$.

The global-local \m\ of $(\rp, \scb, \nu, T)$ w.r.t.\ $\scv$, $\gou(\leb)$ and 
$L^1(\nu)$ follows from Theorem \ref{thm-main-rp} via Lemma 
\ref{lem-primo-ultimo}, applied with $\go = \gou(\leb)$.
\qed

\subsection{Maps of the unit interval}
\label{subs-maps-ui}

We now consider maps $T: \ui \into \ui$ which are Markov for the partition 
$\{ I_j \}_{j=0}^{N-1}$, where $I_j = (a_j, a_{j+1})$ and $0 = a_0 < a_1 
< \ldots < a_{N-1} < a_N = 1$. Once again, we mention now and then no more
that $\{ I_j \}$ is only a partition mod $\leb$, which here denotes the Lebesgue 
\me\ on $\ui$. These are the assumptions on $T$:

\begin{itemize} 
\item[(B1)] $T |_{I_j}$ is a bijective map onto $(0,1)$, and has a $C^2$ 
extension $\br_j: [a_j, a_{j+1}] \into [0,1]$ (this means that $\br_j$ is $C^2$
up to the boundary of its domain).

\item[(B2)] There exists $\Lambda > 1$ such that $| \br'_j | \ge \Lambda$, 
for all $j \ge 1$. Also, $\br'_0(\x) > 1$, for all $\x \in (0,a_1]$.

\item[(B3)] As $\x \to 0^+$, $T'(\x) = 1 + \chi \x^p + o(\x^p)$, for 
some $p \ge 1$ and $\chi>0$.

\item[(B4)] $T$ is exact w.r.t.\ $\leb$ (again, the Lebesgue \me\ on $\ui$).
\end{itemize}

\begin{remark}
  It was proved in \cite{t83} that, assuming (B1)-(B2), $\br_0'(0)=1$ (which 
  is part of (B3)) and 
  \begin{itemize}
  \item[(B4)'] $T'' > 0$ in a neighborhood of 0,
  \end{itemize}
  $T$ is conservative and exact (w.r.t.\ $\leb$). Moreover, assuming (B1)-(B3),
  and (B4)', there exists a unique $T$-invariant \me\ $\mu$ such that 
  $\frac{d\mu}{d\leb}(\x) \sim \x^{-p}$, for $\x \to 0^+$. This comes from 
  \cite[\S1]{t01}, since (B3) implies that 
  \begin{equation} \label{taylor-t}
    T(\x) = \x + \frac \chi {p+1}\, \x^{p+1} + o(\x^{p+1}) ; 
  \end{equation}
  cf.\ also \cite{t80} and \cite[Thm.~2.1]{bl}.
\end{remark}

The above remark shows how the \me\ $\lambda_p$, defined by
\begin{equation} \label{def-lambdap}
  \frac{d\lambda_p}{d\leb}(\x) = \frac 1 {\x^p},
\end{equation}
is especially relevant for maps satisfying (B3). The next theorem, which is the 
main result of this section, reflects this. The case $p=1$ presents certain issues 
that lead us to consider instead (uniformly) global \ob s for the \me\ 
$\lambda_{1+}$, defined by 
\begin{equation} \label{def-lambda1plus}
  \frac{d \lambda_{1+}}{dm}(\x) = -\frac{\log \x}{\x}.
\end{equation}

\begin{theorem} \label{thm-main-ui}
  Let $T$ be a map satisfying (B1)-(B4). If $p > 1$, then $T$ with $\lambda_p$ is 
  fully global-local \m\ w.r.t.\ $\scv = \scv_\mathrm{unif}$. If $p=1$, then $T$ with 
  $\lambda_1$ is global-local \m\ w.r.t.\ $\scv = \scv_\mathrm{unif}$, 
  $\gou(\lambda_{1+})$ and $L^1(\lambda_1)$.
\end{theorem}

\begin{remark}
  The purpose of $\lambda_{1+}$ is simply to have a \me\ whose singularity at
  $0$ is slightly stronger than $1/\x$. The specific form (\ref{def-lambda1plus}) 
  is unimportant, except that it behaves well with the many technical arguments 
  of the proof of Theorem \ref{thm-main-ui} (Section \ref{subs-pf-main-ui}). Notice
  that the case $p=1$ was also special for the main theorem of \cite{bl} (Thm.\ 2.4) 
  about general global \ob s.
\end{remark}

Our study of the restrictions of infinite-volume averages shows that, if $T$ is 
fully global-local \m\ for $\lambda_p$, $p>1$, then it is global-local \m, though 
not fully, for many other measures, including $\lambda_q$, for $1 \le q < p$.

\begin{corollary} \label{cor-main-ui}
  For $1 \le q < p$, let $\nu$ be an infinite, locally finite \me, absolutely continuous
  w.r.t.\ the Lebesgue \me\ $\leb$, with $\frac{d\nu}{d\leb}(\x) 
  \sim \x^{-q}$, as $\x \to 0^+$. Then $\overline{\lambda_p}$ is a strict restriction 
  of $\overline{\nu}$. Furthermore, under the same hypotheses as in Theorem 
  \ref{thm-main-ui}, $(\ui, \scb, \nu, T)$ is global-local \m\ w.r.t.\ $\scv$, 
  $\gou(\lambda_p)$ and $L^1(\nu)$.
\end{corollary}

\proof Using Proposition \ref{prop:restriction} and Remark \ref{rm:simpl-cond}, it is 
easy to see that $\overline{\lambda_p}$ is a restriction of $\overline{\lambda_q}$, 
and by Proposition \ref{prop2} the same holds for $\overline \nu$. To show that 
the restriction is strict, it suffices to produce a global observable in 
$\gou(\lambda_q) \setminus \gou(\lambda_p)$. For $q > 1$, we may take
the \fn\ $F_q(\x) := \cos(\x^{-q+1})$, because, for all $[a,b] \subset \ui$, 
\begin{equation}
  \frac1 {\lambda_q([a,b])} \int_a^b F_q \, d\lambda_q = \frac{\sin (a^{-q}) - 
  \sin (b^{-q})} {\lambda_q([a,b])\, q},
\end{equation}
implying that $F_q \in \gou(\lambda_q)$, with $\overline{\lambda_q} (F_q) = 0$. 
On the other hand, let $r > 0$ and $V_b = [b-c,b] \subset \ui$, such that 
$\lambda_p (V_b) = r$. A simple calculation shows that 
$c = b - (b^{1-p} + r(p-1) )^{1/(1-p)}$. In particular, $c \sim r b^p$, for $b \to 0^+$. 
Therefore, since $q < p$, the function $F_q$ is approximately constant on $V_b$, 
when $b$ is very close to $0$, giving 
\begin{equation}
  \underset{b \to 0^+}{\limsup / \liminf} \ \frac1r \int_{V_b} F_q \, d\lambda_p = \pm 1.
\end{equation}

For $q=1$, a similar argument holds for the global observable $F_1(\x) := 
\cos( \log(\x^{-1}))$, yielding that $\gou(\lambda_1) \setminus \gou(\lambda_p)$ is 
non-empty for all $p>1$.

The global-local \m\ of $(\ui, \scb, \nu, T)$ as in the statement of the corollary
follows once again from Lemma \ref{lem-primo-ultimo} (with $\go = 
\gou(\lambda_p)$).
\qed

\section{Proofs}
\label{sec-proofs}

In this section we prove our main results and Proposition \ref{prop:restriction}, that
was left behind. We start with the latter.

\subsection{Proof of Proposition \ref{prop:restriction}}

Let $F \in \gou(\mu)$ and, without loss of generality, assume that $\avg(F) = 0$. 
We want to prove that $\overline{\nu}(F)$ exists and equals 0.

Let us start by noting that, for any $V \in \scv$,
\begin{equation}
\begin{split}
  & \left| \int_{V} F \frac{d\nu}{d\mu} \, d\mu - \frac{\nu(V)}{\mu(V)} \int_{V}F \, d\mu
  \right| = \left|\int_{V} F \left( \frac{d\nu}{d\mu} - \frac{\nu(V)}{\mu(V)} \right) d\mu
  \right| \\
  &\qquad \le \|F\|_{\infty} \left(\sup_{V} \frac{d\nu}{d\mu} - \inf_{V} \frac{d\nu}{d\mu} 
  \right) \mu(V) \\
  &\qquad = \|F\|_{\infty} \left[ \sup_{V} \frac{d\nu}{d\mu} \left(\inf_{V} \frac{d\nu}{d\mu}
  \right)^{-1} - 1 \right] \mu(V) \, \inf_{V} \frac{d\nu}{d\mu},
\end{split}
\end{equation}
where we have used that, for any $x \in V$, 
\begin{equation}
  \left| \frac{d\nu}{d\mu}(x) -\frac{\nu(V)}{\mu(V)} \right| = \left| \frac{d\nu}{d\mu}(x) -
  \frac1 {\mu(V)} \int_V \frac{d\nu}{d\mu} \, d\mu \right| \le \sup_V \frac{d\nu}{d\mu} - 
  \inf_V \frac{d\nu}{d\mu} .
\end{equation}
Dividing both sides by $\nu(V)$ and using that $\mu(V)\inf_{V}\frac{d\nu}{d\mu} \le 
\nu(V)$, we have the following inequality that we will use later,
\begin{equation} \label{eq:resp-flat}
  \left| \frac1 {\nu(V)} \int_{V} F\, d\nu - \frac1 {\mu(V)} \int_{V} F\, d\mu \right| \le 
  \|F\|_{\infty} \left[\sup_{V} \frac{d\nu}{d\mu} \left(\inf_{V} \frac{d\nu}{d\mu} \right)^{-1}
  -1\right].
\end{equation}

Fix $\eps>0$. Since $\avg(F) = 0$, there exists $r=r(\eps)>0$ so large that
\begin{equation} \label{eq:mu-avg}
  \sup_{V \in \scv_\mu(r)}\left| \frac1r \int_V F \, d\mu \right| \le \frac{\eps}{4}
\end{equation}
and (\ref{eq:flat})-(\ref{eq:comparability}) hold. So, by (\ref{eq:flat}), there exists 
$M \in \ps$ such that, for all $V \in \scv_{\mu}(r)$, $V \subset U_M$,
\begin{equation} \label{eq:small-derivative}
 \left[ \sup_{V} \frac{d\nu}{d\mu} \left(\inf_{V} \frac{d\nu}{d\mu} \right)^{-1} - 1 
 \right] \le \frac{\eps}{4\|F\|_{\infty}}.
\end{equation}
For $R > 0$, which one should think of as much larger than $r$, let us consider 
$W \in \scv_\nu(R)$. Partition $W$ into adjacent intervals $V_0, V_1, \ldots, V_l$,
from left to right in the case $\ps = \rp$, and from right to left in the case $\ps = 
\ui$, such that $\mu(V_i)=r$, for $i=1, \ldots, l$, and $\mu(V_0) \le r$; see 
Fig.~\ref{fig1}. Let $l'$ be the smallest integer between 1 and $l$ such that $V_{l'} 
\subset U_M$. Of course, $V_i \subset U_M$ for all $i\ge l'$.

\newfig{fig1}{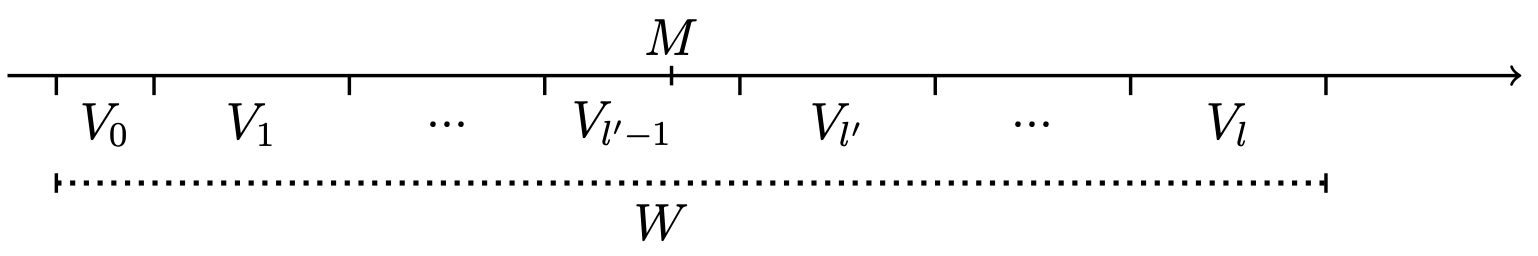}{12cm}{Partition of $W$ for the proof of Proposition 
\ref{prop:restriction} (case $\ps = \rp$).}

Now write
\begin{equation} \label{eq:subdivision}
  \int_W F \, d\nu = \sum_{i=0}^{l'-1} \int_{V_i} F \, d\nu + \sum_{i=l'}^l 
  \int_{V_i} F \, d\nu.
\end{equation}
For $i\ge l'$, since $V_i \subset U_M$,
\begin{equation} \label{eq:good-terms}
  \left| \frac1 {\nu(V_i)} \int_{V_i} F\, d\nu - \frac1 {\mu(V_i)} \int_{V_i} F\, d\mu
  \right| \le \frac{\eps}{4},
\end{equation}
where we have used (\ref{eq:resp-flat}) and (\ref{eq:small-derivative}).
The estimate (\ref{eq:good-terms}), together with (\ref{eq:mu-avg}) and the fact 
that $V_i \in \scv_\mu(r)$, implies that, for all $i \ge l'$,
\begin{equation} \label{eq:almost-done}
  \left |\int_{V_i} F \, d\nu \right| \le \nu(V_i) \, \frac{\eps}{2}.
\end{equation}
As for the other terms in the r.h.s.\ of (\ref{eq:subdivision}), observe that
\begin{equation} \label{eq:bad-intervals}
  \sum_{i=0}^{l'-1} \int_{V_i} F \, d\nu \le \|F\|_{\infty} \left( \nu(\ps \setminus U_M) 
  + \sup_{V \in \scv_\mu(r)} \nu(V) \right),
\end{equation}
having estimated $| \int_{V_{l'-1}} F d\nu | \le \|F\|_{\infty} 
\sup_{V \in \scv_\mu(r)} \nu(V)$. Notice that this estimate, and thus 
(\ref{eq:bad-intervals}), also holds in the case $l'=1$, because $\mu(V_0) \le r$.
Setting
\begin{equation}
  \overline R = \frac2 \eps \, \|F\|_{\infty} \left(\nu(\ps \setminus U_M) + 
  \sup_{V \in \scv_\mu(r)} \nu(V)\right)
\end{equation}
one has that for any $R \ge \overline R$, the l.h.s.\ of (\ref{eq:bad-intervals}) 
does not exceed $R\eps/2$. This fact, (\ref{eq:subdivision}) and 
(\ref{eq:bad-intervals}) imply that, for all $R \ge \overline R$ and 
$W \in \scv_\nu(R)$,
\begin{equation}
  \left| \int_{W} F\, d\nu\right| \le \frac{R\eps}{2} + \frac{\eps}{2} \sum_{i=l'}^l
  \nu(V_i) \le R \eps,
\end{equation}
which is what was to be proved.
\qed

\subsection{Proof of Theorem \ref{thm-main-rp}}

Here and for the rest of the paper we denote by $\ibr_j := \br_j^{-1}$ the 
\emph{inverse branches} of $T$, for both cases $\ps = \rp$ and $\ps = \ui$.

We begin the proof of Theorem \ref{thm-main-rp} with a number of simple 
lemmas on the properties of the branches and inverse branches of $T$.

\begin{lemma} \label{lem:dist-0}
  For $j \ge 1$, either $\ds \lim_{x \to a_j^-} \tau'_j (x) = +\infty$ or 
  $\ds \lim_{x \to a_{j+1}^+}  \tau'_j (x) = -\infty$, depending on $\br_j$ being 
  increasing or decreasing, respectively.
\end{lemma}

\proof Let us denote with $x \to \bullet$ either limit $x \to a_j^-$ or 
$x \to a_{j+1}^+$, depending on the two cases $\br_j' > 0$ or $\br_j' < 0$, as in
the statement of the lemma. By (A4),
\begin{equation}
  \lim_{x \to \bullet} \left(\frac{1}{\br'_j(x)}\right)' = 0.
\end{equation}
Therefore $1/\br_j'(x)$ has a limit, for $x \to \bullet$, in $\R$. Then, so does 
$\br_j'(x)$, in $\R \cup \{ \pm \infty \}$ (because $\br_j'$ never changes sign). The 
fact that $\br_j$ is monotonic and surjective onto $\rp$ forces the conclusion of the 
lemma.
\qed

\begin{lemma} \label{lem:distortion}
  For all $j$,
  \begin{equation} \label{eq:sup}
     \sup_\rp \, \frac{ |\ibr''_j| } { |\ibr'_j|(1-|\ibr'_j|) } < \infty.
  \end{equation}
  Furthermore,
  \begin{equation} \label{eq:limit-zero}
     \lim_{x \to +\infty} \, \frac{\ibr''_0(x)} {\ibr'_0(x)(1-\ibr'_0(x))} = 0
  \end{equation}
  and, for all $j\ge 1$,
  \begin{equation} \label{eq:limit-j}
     \lim_{x \to +\infty} \ibr'_j(x) = \lim_{x \to +\infty} \, \frac{\ibr''_j(x)}{\ibr'_j(x)} = 0.
  \end{equation}
\end{lemma}

\proof Since $\ibr'_j = \frac{1}{\br'_j}\circ \ibr_j$ and $\ibr''_j = -\frac{\br''_j}
{(\br'_j)^3} \circ \ibr_j$, one gets
\begin{equation} \label{eq:to-bound}
   \frac{ |\ibr''_j| }  {|\ibr'_j| \, (1-|\ibr'_j|) } = \frac{ |\br''_j| } { |\br'_j| \, (|\br'_j|-1) }
   \circ \ibr_j. 
\end{equation}

Let us prove (\ref{eq:sup}) for $j=0$. From (A1)-(A3), the 
\fn\ $\frac{|\br_0''|}{\br_0' (\br_0'-1)}$ is continuous on $[a_0, +\infty)$ (observe
that $\br_0'>1$ up to and including $a_0$) and vanishes at $+\infty$, implying 
that (\ref{eq:to-bound}) is bounded for $j=0$. In addition, (A3) gives 
(\ref{eq:limit-zero}).

For the case $j\ge 1$, let us rewrite the r.h.s.\ of (\ref{eq:to-bound}) as
\begin{equation} \label{eq:to-bound2}
   \frac{ |\br''_j| } { |\br'_j|^2 \left( \ds 1-\frac1 {|\br'_j|} \right)} \circ \ibr_j.
\end{equation}
Since $|\br'_j| \ge \Lambda > 1$, a sufficient condition for (\ref{eq:sup}) is an 
upper bound for $|\br''_j| / |\br'_j|^2$, which we readily get from (A4) 
and the continuity of $\br_j''$. Finally, since
\begin{equation} 
   \frac{\ibr''_j}{\ibr'_j} = - \frac{\br''_j}{ (\br'_j)^2} \circ \ibr_j, 
\end{equation}
(\ref{eq:limit-j}) follows from Lemma \ref{lem:dist-0} and (A4).
\qed

\begin{lemma} \label{lem-br0-x}
  For $x \in I_0$, $\br_0(x) < x$.
\end{lemma}

\proof If the claim is false, there must exist $a \in I_0$ such that $\br_0(a) = a$.
Since $\br_0' > 1$, it follows that $\br_0(x) > x$, for all $x>a$. Thus, chosen some 
$x_0 > a$ and denoted $x_1 := T(x_0)$, we have that $x_1>x_0$ and
$\{ T^n [x_0, x_1) \}_{n\in \N}$ is a partition of $[x_0, +\infty)$. This contradicts
the exactness of $T$ via the so-called Miernowski-Nogueira criterion \cite{mn}.
(See Appendix A.1 of \cite{lsimple} for a generalization of said criterion. In the 
language of that reference, $[x_0, x_1)$ is not asymptotically intersecting and 
therefore $T$ cannot be exact.)
\qed

Let us introduce the main objects that will be used in the proof of Theorem
\ref{thm-main-rp}:

\paragraph{Cones of local \ob s:} For $M, D, \eps > 0$, set
\begin{equation} \label{cones-rp}
   \mathcal C_{M,D,\eps} := \rset{g \in L^1(\rp, \leb) \cap C^1(\rp)} {g>0, 
   \frac{|g'|}{g} \le D \mbox{ and } \frac{|g'(x)|}{g(x)} \le \eps, \forall x \ge M}.
\end{equation}
This is a \emph{cone}, in the sense that $g \in \mathcal C_{M,D,\eps}$ and 
$c>0$ implies $cg \in \mathcal C_{M,D,\eps}$.

\paragraph{Transfer operator:} We define $\trop: L^1(\leb) \into L^1(\leb)$ via the 
identity
\begin{equation} 
  \forall F \in L^{\infty}(\leb), g \in L^1(\leb), \qquad \int_\rp (F \circ T) g \, d\leb 
  = \int_\rp F \, \trop g \, d\leb.
\end{equation}
In this case, where the reference \me\ is the Lebesgue \me, $\trop$ is also 
referred to as the \emph{Perron-Frobenius operator} and its explicit formula is
notoriously:
\begin{equation} \label{trop-expl}
  \trop g  = \sum_{j=0}^{N-1} |\ibr'_j| \, (g\circ \ibr_j) = \sum_{j=0}^{N-1} \sigma_j 
  \ibr'_j \, (g\circ \ibr_j),
\end{equation}
where $\sigma_j := \mathrm{sgn}(\ibr_j')$ (recall that $\ibr_j$ is monotonic).

Now set
\begin{equation} \label{d-under}
  \underline{D} := \max_{0 \le j \le N-1} \, \sup_\rp \, \frac{ |\ibr''_j| }
  { |\ibr'_j| \, (1- |\ibr'_j|) } ,
\end{equation}
which is finite by Lemma \ref{lem:distortion}. The following result is a 
crucial statement that says that the dynamics preserves densities that are 
``almost flat at infinity''.

\begin{lemma} \label{lem:cone-invariance}
  Let $D \ge \underline{D}$. For all $\eps>0$, there exists $M>0$ such 
  that $\trop \mathcal C_{M, D, \eps} \subseteq \mathcal C_{M, D, \eps}$.
\end{lemma}

\proof Let $g \in C^1$, $g > 0$. We first establish the inequality 
\begin{equation} \label{eq:deriv-bound}
   |(\trop g)' | \le \max_j \left( \left|\frac{\ibr''_j} {\ibr'_j} \right| +
   \frac{| (g' \circ \ibr_j) \ibr'_j |} {g \circ \ibr_j} \right) \trop g,
\end{equation}
that will be used repeatedly in the present proof. A simple computation out of
(\ref{trop-expl}) yields
\begin{equation} \label{eq:deriv-t-op}
   (\trop g)' = \sum_{j=0}^{N-1} \left( \sigma_j \ibr''_j \, (g \circ \ibr_j) + 
   \sigma_j (\ibr'_j)^2 \, (g'\circ \ibr_j) \right).
\end{equation}
If in each of the above summands we divide and multiply the left term by 
$\ibr'_j$ and the right term by $g \circ \ibr_j$, we obtain (\ref{eq:deriv-bound})
after a triangular inequality.

Now let $\eps > 0$. By (\ref{eq:limit-zero}) and (\ref{eq:limit-j}) there exists $M>0$
such that, for all $x \ge M$,
\begin{align}
  \label{eps-est1}
  \frac{ |\ibr''_0(x)| } {\ibr'_0(x)(1-\ibr'_0(x))} &\le \eps; \\[3pt]
  \label{eps-est2}
  \max_{j \ge 1}\left( \left| \frac{\ibr''_j}{\ibr'_j}\right| +D |\ibr'_j| \right) &\le \eps.
\end{align}
Let us show that $\mathcal C_{D, M, \eps}$ is $\trop$-invariant, as in the 
claim of Lemma \ref{lem:cone-invariance}. Take $g \in \mathcal C_{M, D, \eps}$. 
Since $\ibr_j \in C^2$, $\trop g \in C^1$. Also, the surjectivity of $T$ and $g>0$ 
imply $\trop g>0$. Moreover, by (\ref{eq:deriv-bound}), $|g'|/g \le D$ and
$D \ge \underline{D}$, cf.\ (\ref{d-under}), we obtain
\begin{equation} 
   |(\trop g)' | \le \max_j \left( \left| \frac{\ibr''_j}{\ibr'_j} \right| + D |\ibr'_j| \right) 
   \trop g \le D \, \trop g.
\end{equation}

At this point, take $x \ge M$. By Lemma \ref{lem-br0-x}, $\ibr_0(x) > x \ge M$.
The properties of $g \in \mathcal C_{D, M, \eps}$ imply that
\begin{equation} 
\begin{split}
   |(\trop g)'(x)| &\le \max_j \left( \left| \frac{\ibr''_j(x)} {\ibr'_j(x)} \right| +\left|
   \frac{g' (\ibr_j(x)) \, \ibr'_j(x)} {g(\ibr_j(x))} \right| \right) \trop g(x)  \\
   &\le \max \left\{ \left( \frac{ |\ibr''_0(x)| } {\ibr'_0(x)} + \eps \ibr'_0(x) 
   \right), \, \max_{j \ge 1} \left( \left| \frac{\ibr''_j(x)} {\ibr'_j(x)} \right| +D 
   |\ibr'_j(x)| \right) \right\} \trop g(x) \\
  &\le \eps \, \trop g(x) ,  
\end{split}
\end{equation}
where the last inequality comes from (\ref{eps-est1})-(\ref{eps-est2}). This
concludes the proof of Lemma \ref{lem:cone-invariance}.
\qed

The next lemma formalizes the intuitive fact that if a local \ob\ is almost flat
in a neighborhood of infinity, its integral against a uniformly global \ob\ in that 
neighborhood results approximately in the infinite-volume average of the
global \ob.

\begin{lemma} \label{lem:control-avg}
  Let $F \in \gou(\leb)$, $\avgleb(F) = 0$. For each $\delta>0$, there 
  exists $\eps > 0$ such that for, every $g \in C^1(\rp)$ and $M\ge 0$ with
  \begin{itemize}
  \item[(i)] $g > 0$;
  \item[(ii)] $\leb(g) = 1$;
  \item[(iii)] $\ds \frac{|g'(x)|}{g(x)} \le \eps$, for all $x \ge M$,
  \end{itemize}
  one has
  \begin{displaymath}
     \left |\int_M^\infty F g \, d\leb \right| \le \delta.
  \end{displaymath}
\end{lemma}

\proof Let us assume that $F$ is not almost everywhere $0$, otherwise the
statement is trivial. By the hypothesis $\avgleb(F) = 0$, given $\delta>0$
there exists $r>0$ such that 
\begin{equation} \label{eq:globavg-lem}
   \sup_{V \in \scv_\leb(r) } \left| \frac1r \int_V F \, d\leb \right| \le \frac \delta 2. 
\end{equation}
Set $\eps := \frac \delta {2r \| F\|_\infty}$. 
For any $g$ and $M$ verifying \emph{(i)-(iii)}, we decompose
\begin{equation} \label{eq:split-integral}
  \int_M^\infty F g \, d\leb = \sum_{k=0}^\infty \int_{M + kr}^{M + (k+1)r} F g \, 
  d\leb.
\end{equation}
For $k \in \N$, call $x_{\mathrm{min},k}$ a minimum point of $g$ in 
$[M + kr, M + (k+1)r]$ (recall that $g$ is continuous). By (\ref{eq:globavg-lem}), 
\emph{(i)}, \emph{(iii)} and the definition of $\eps$,
\begin{equation} \label{eq:small-pieces}
\begin{split}
  \left| \int_{M + kr}^{M + (k+1)r} F g \, d\leb \right| &= 
  \left| \int_{M + kr}^{M+ (k+1)r} F(x) \left( g(x_{\mathrm{min},k}) + 
  \int_{x_{\mathrm{min},k}}^x g'(t) \, dt \right) dx \right| \\
  &\le \frac \delta 2 r g(x_{\mathrm{min},k}) + r \|F\|_\infty \, \eps 
  \left(\int_{M + kr}^{M+(k+1)r} g(t) \, dt \right) \\
  &\le \left(\frac \delta 2 + \frac \delta 2 \right) \left(\int_{M + rK}^{M+(k+1)r} g(t) \, dt 
  \right).
\end{split}
\end{equation}
We obtain the assertion of Lemma \ref{lem:control-avg} by 
(\ref{eq:split-integral})-(\ref{eq:small-pieces}), using \emph{(ii)}.
\qed

We are now ready for the final arguments of the proof of Theorem 
\ref{thm-main-rp}. It suffices to establish that 
\begin{equation} \label{goal1}
  \lim_{n \to \infty} \, \leb((F \circ T^n) g) = 0,
\end{equation}
for all $F \in \gou(\leb)$ with $\avgleb(F) = 0$ and all $g \in L^1(\leb)$. 
If $\leb(g)=0$, the above follows from the weak form of global-local \m\ called
(GLM1), which is a simple consequence of the exactness of $T$ 
\cite[Thm.~3.5]{lpmu}.%
\footnote{This theorem and Lemma 3.6 in the same reference, which we use 
momentarily, were stated for measure-preserving systems, with another 
assumptions that ensured that the infinite-volume average was 
dynamics-invariant. None of these assumptions are needed for the proofs of 
the statements invoked here. See also Appendix A.1 of \cite{bl}.} 
So we turn to the case $\leb(g) \ne 0$. By the linearity of (\ref{goal1}) in $g$, we 
may assume $\leb(g)=1$. Since $T$ is exact, it is enough to show that, for 
all $\delta>0$, there exists $g_\delta \in L^1(\leb)$, which may depend on $F$, 
such that $\leb(g_\delta)=1$ and
\begin{equation} \label{goal2}
  \limsup_{n \to \infty} \, | \leb((F \circ T^n) g_\delta) | \le 2\delta.
\end{equation}
In fact, \cite[Lem.~3.6]{lpmu} states that if the above conditions holds for $g_\delta$, 
then it holds for \emph{all} $g$ with $\leb(g)=1$. The arbitrariness of $\delta$ then
implies (\ref{goal1}) for all such $g$.

It remains to prove (\ref{goal2}). Fix $\delta>0$ and let $\eps>0$ be given by 
Lemma \ref{lem:control-avg} for the chosen value of $\delta$. Then let $M>0$ 
be given by Lemma \ref{lem:cone-invariance} for the selected value of $\eps$
and any fixed value of $D \ge \underline{D}$. Let $g_\delta$ be any element
of $\mathcal C_{M,D,\eps}$ with $\leb(g_\delta) = \| g_\delta \|_1 =1$ (such an 
element evidently exists). Then, for all $n \in \N$,
\begin{equation}
\begin{split}
  | \leb((F \circ T^n) g_\delta) | &\le | \leb( (1_{[0,M]} F) \circ T^n) g_\delta) | + 
  | \leb( (1_{(M, \infty)} F) \circ T^n) g_\delta) | \\
  &= | \leb( (1_{[0,M]} F) \circ T^n) g_\delta) | + 
  \left| \int_M^\infty F (\trop^n g_\delta) \, dm \right| .
\end{split}
\end{equation}
We conclude the proof of Theorem \ref{thm-main-rp} if we are able to bound
each the above terms with $\delta$, for $n$ large. The first bound comes from 
the exactness of $T$, which implies full \emph{local-local \m} (LLM) w.r.t.\ $\leb$, 
namely,
\begin{equation} \label{eq:exactness-1}
  \forall f \in L^\infty(\leb) \cap L^1(\leb), \, \forall g \in L^1(\leb), \quad 
  \lim_{n \to \infty} \, \leb((f \circ T^n) g) = 0;
\end{equation}
cf.\ \cite[Thm.~3.5]{lpmu}. The second bound holds for all $n$ and comes directly 
from Lemma \ref{lem:control-avg}, since $\trop^n g_\delta \in 
\mathcal C_{M,D,\eps}$  (by Lemma \ref{lem:cone-invariance}) and 
$\leb( \trop^n g_\delta  ) =  \| \trop^n g_\delta \|_1 =  \| g_\delta \|_1 = 1$. 
\qed

\subsection{Proof of Theorem \ref{thm-main-ui}}
\label{subs-pf-main-ui}

Once again, we start with a few calculus lemmas.

\begin{lemma} \label{lem:second-derivative}
  For $p>1$, $\ibr_0''(\x) \to 0 $ as $\x \to 0^+$.
\end{lemma}

\proof Since $\ibr_0'' = -\frac{\br_0''}{(\br_0')^3} \circ \ibr_0$, and 
$\lim_{\x \to 0^+}\br_0'(\x) =1$, it is enough to prove that 
$\lim_{\x \to 0^+} \br_0''(\x) = 0$. By (B1), $\br_0''(\x) \to c \in \R$, as $\x \to 0^+$. 
If $c \ne 0$, the Taylor expansion of $\br_0'$ at 0 would be $\br_0' = 1 +c\x + 
o(\x)$, contradicting (B3).
\qed

Let $\ell_p, \ell_{1+}: (0,1] \into \R^+$ denote the (infinite) densities of the \me s
$\lambda_p, \lambda_{1+}$, respectively:
\begin{align}
  \ell_p(\x) &:= \frac 1 {\x^p}; \label{ellp} \\
  \ell_{1+}(\x) &:= -\frac{\log \x}{\x}.
\end{align}
As already mentioned in Section \ref{subs-maps-ui}, we will work with the \fn\
$\ell_p$, for $p>1$, and $\ell_{1+}$, for $p=1$. Since we prove
the same statements for all these \fn s, it is convenient to introduce the
common notation%
\footnote{The reader might think that a better notation would be to define $\ell_p$
as in (\ref{ellp}) and $\ell_1(\x) := -\log \x/\x$, but unfortunately we also need 
the function $\ell_1(\x) := 1/\x$.}
\begin{equation} 
  p\symb = \begin{cases}
  p, & \text{if } p >1; \\
  1+, & \text{if } p=1.
\end{cases}
\end{equation}
The analogue of Lemma \ref{lem:distortion} for interval maps is the combination of
the upcoming Lemmas \ref{lem:dist-int-2} and \ref{lem:dist-int-3}, which require 
a few preparatory results.

\begin{lemma} \label{lem:distortion-int}
  Under the above assumptions,
  \begin{displaymath}
    \lim_{\x \to 0^+} \frac{ 1-\ibr'_0(\x) \, \ds \frac{ \ell_{p\symb}(\x) \, 
    \ell'_{p\symb}(\ibr_0(\x)) } { \ell'_{p\symb}(\x) \, \ell_{p\symb}(\ibr_0(\x)) } } 
    {1 - \ibr'_0(\x) \, \ds \frac{ \ell_{p\symb} (\ibr_0(\x))} {\ell_{p\symb}(\x) } } = p.
  \end{displaymath}
\end{lemma}

\proof Set $\kappa := \frac{\chi}{p+1}$. By (B3) and (\ref{taylor-t}), we have, as 
$\x \to 0^+$,
\begin{align} 
  \ibr_0' (\x) &= 1 - \chi \x^p + o(\x^p) ; \\
  \ibr_0(\x) &= \x - \kappa \x^{p+1} + o(\x^{p+1}). \label{eq:asimp-first-inv}
\end{align}

Let us first consider the case $p>1$. The above expansions give
\begin{equation} 
\begin{split}
  \frac{ 1-\ibr'_0(\x) \, \ds \frac{ \ell_{p}(\x) \, 
  \ell'_{p}(\ibr_0(\x)) } { \ell'_{p}(\x) \, \ell_{p}(\ibr_0(\x)) } } 
  {1 - \ibr'_0(\x) \, \ds \frac{ \ell_{p} (\ibr_0(\x))} {\ell_{p}(\x) } } &=
  \frac{1- (1- \chi \x^p + o(\x^p)) \, \ds \frac{\x}{\ibr_0(\x)}} {1- (1-\chi\x^p + o(\x^p))
  \, \ds \frac{\x^p}{ (\ibr_0(\x))^p }} \\
  &=  \frac{1- (1- \chi \x^p + o(\x^p) ) \, (1-\kappa\x^p + o(\x^p))^{-1}}
  {1- (1-\chi\x^p + o(\x^p)) \, (1-\kappa\x^p + o(\x^p))^{-p}} \\
  &= \frac{\chi - \kappa + o(1) } {\chi - p \kappa + o(1) } \\
  &=p +o(1),
\end{split}
\end{equation}
as claimed. In the case $p=1$ we have instead
\begin{equation} 
  \frac{ \ell_{1+}(\x) \, \ell'_{1+}(\ibr_0(\x)) } { \ell'_{1+}(\x) \, \ell_{1+}(\ibr_0(\x)) } 
  = \frac{\x}{\ibr_0(\x)} \left[ \frac{\log \x} {\log \ibr_0 (\x)} \, 
  \frac{\log \ibr_0(\x) -1}{\log \x -1} \right].
\end{equation}
Therefore, since $\chi - \kappa \ne 0$,
\begin{equation} \label{eq:p=1}
\begin{split}
  & \frac{ 1-\ibr'_0(\x) \, \ds \frac{ \ell_{1+}(\x) \, 
  \ell'_{1+}(\ibr_0(\x)) } { \ell'_{1+}(\x) \, \ell_{1+}(\ibr_0(\x)) } } 
  {1 - \ibr'_0(\x) \, \ds \frac{ \ell_{1+} (\ibr_0(\x))} {\ell_{1+}(\x) } } \\
  &\qquad = \frac{1- (1- \chi \x + o(\x)) \, \ds \frac{\x}{\ibr_0(\x)} \left[ \frac{\log \x} 
  {\log \ibr_0 (\x)}  \, \frac{\log \ibr_0(\x) -1}{\log \x -1} \right] } {1- (1-\chi\x + o(\x))
  \, \ds \frac{\x}{\ibr_0(\x)} \, \frac{\log \ibr_0(\x)} {\log \x} } \\
  &\qquad = \frac{1- (1- \chi \x + o(\x)) \, (1+ \kappa\x + o(\x)) \ds 
  \left[ \frac{\log \x} {\log \ibr_0 (\x)} \, \frac{\log \ibr_0(\x) -1}{\log \x -1} \right] } 
  {1- (1-\chi\x + o(\x)) \, (1+ \kappa\x + o(\x)) \, \ds \frac{\log \ibr_0(\x)} {\log \x} }. 
\end{split}
\end{equation}
Using (\ref{eq:asimp-first-inv}) with $p=1$, we obtain
\begin{align}
  \frac{ \log \ibr_0(\x)} {\log \x} &= 1 + o(\x); \label{eq:fancy-orders} \\
  \frac{\log \ibr_0(\x) -1}{\log \x -1} &= 1 + o(\x),
\end{align}
so that (\ref{eq:p=1}) is equal to $1 + o(1)$, proving Lemma 
\ref{lem:distortion-int}.
\qed

The next lemma is obtained through (some of) the same estimates as 
presented above, so we omit its proof.

\begin{lemma} \label{lem:another-dist}
Under the above assumptions,
\begin{displaymath}
  \lim_{\x \to 0^+} \left[ \ell_{p\symb} (\x) \left( 1 - \ibr'_0(\x) \, \ds 
  \frac{\ell_{p\symb} (\ibr_0 (\x))} {\ell_{p\symb} (\x)} \right) \right] =
  \begin{cases}
    \ds \frac \chi {p+1}, & \text{if } p >1; \\[8pt]
    +\infty, & \text{if } p=1.
  \end{cases}
\end{displaymath}
\end{lemma}

It will be convenient below to use modifications of the densities $\ell_{p\symb}$ 
defined as follows:
\begin{equation} 
  \label{eq:modification-def}
  z_{p\symb} (\x) := 
  \begin{cases}
    \ell_{p\symb}(\x), & \x \in (0, \eta/2]; \\
    \vartheta(\x), & \x \in (\eta/2, \eta); \\
    c & \x \in [\eta,1],
  \end{cases}
\end{equation}
where $\eta \in (0,a_1)$ will be fixed in Lemma \ref{lem:finite-modification}, $c \in 
(\ell_{p\symb} (\eta), \ell_{p\symb}(\eta/2))$, and $\vartheta_{p\symb} \in 
C^1((\eta/2, \eta))$ is a decreasing \fn\ with
\begin{align}
  \lim_{\x \to (\eta/2)^+} \vartheta(\x) = \ell_{p\symb}(\eta/2), &\qquad 
  \lim_{\x \to \eta^-} \vartheta(\x) = c ; \\
  \lim_{\x \to (\eta/2)^+} \vartheta'(\x) = \ell_{p\symb}'(\eta/2), &\qquad 
  \lim_{\x \to \eta^-} \vartheta'(\x) = 0;\\
  |\vartheta'(\x)| \le |\ell_{p\symb}'(\x)|, &\qquad \forall \x \in (\eta/2, \eta);
\end{align}
see Fig.~\ref{fig2}. Here $\eta, c, \vartheta$ depend on $p\symb$. Observe that, 
by construction,
\begin{equation} \label{eq:mod-properties}
  \ell_{p\symb} \le z_{p\symb} , \qquad  
  |z'_{p\symb}| \le |\ell'_{p\symb}| .
\end{equation}
We also denote by $\zeta_{p\symb}$ the (infinite) \me\ on $\ui$ with density 
$z_{p\symb}$. The \emph{raison d'\^etre} of $z_{p\symb}$ is the following 
property.

\newfig{fig2}{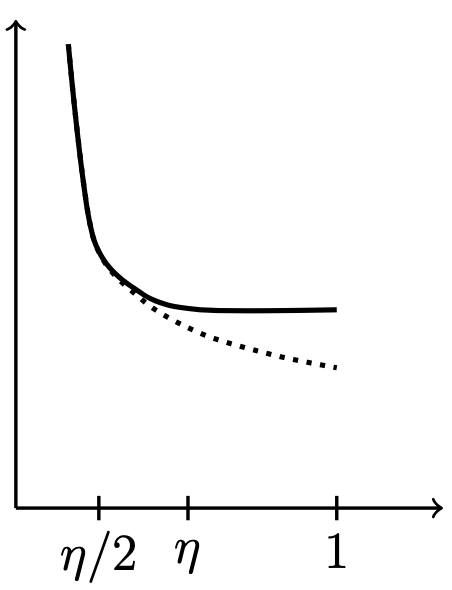}{3.5cm}{Graph of $z_{p\symb}$ (solid line) vs graph of 
$\ell_{p\symb}$ (dashed line)}

\begin{lemma} \label{lem:finite-modification}
  For $\eta$ small enough (depending on $p\symb$) and $j \in \{0, 1, \ldots, N-1\}$,
  \begin{displaymath}
    |\ibr_j'| \, \frac{z_{p\symb} \circ \ibr_j}{z_{p\symb}} < 1.
  \end{displaymath}
  Furthermore, for $j \in \{1, \ldots, N-1\}$, the above l.h.s.\ is less than or equal to 
  $\Lambda^{-1}$.
\end{lemma}

\proof First, observe that, for all $j \ge 1$ and $\x \in \ui$, $\ibr_j(\x) \ge a_1$. By the
monotonicity of $z_{p\symb}$ and $\eta < a_1$,
\begin{equation} 
  z_{p\symb} (\ibr_j(\x)) \le z_{p\symb} (a_1) = c = \min_\ui z_{p\symb}, 
\end{equation}
which, together with $|\ibr_j' |\le \Lambda^{-1}$, proves the second statement
of the lemma (and thus the first, limited to $j\ge 1$).

For the case $j = 0$, notice that, by Lemma \ref{lem:another-dist}, there exists 
$\delta > 0$ such that, for all $\x \in (0,\delta)$,
\begin{equation} \label{eq:mannaggia}
  1 - \ibr'_0(\x) \, \ds \frac{\ell_{p\symb} (\ibr_0 (\x))} {\ell_{p\symb} (\x)} > 0, \quad 
  \text{i.e.,} \quad \ibr'_0(\x) \, \ds \frac{\ell_{p\symb} (\ibr_0 (\x))} {\ell_{p\symb} (\x)} 
  < 1.
\end{equation}
Now choose any $\eta \le \ibr_0(\delta)$. For $\x \in [\delta,1]$, using the fact
that $\ibr_0$ is increasing and $z_{p\symb}$ is decreasing, we have
\begin{equation} \label{eq:finefine}
  \frac{z_{p\symb} (\ibr_0 (\x))} {z_{p\symb} (\x)} \le \frac{z_{p\symb} (\eta)}
  {\min z_{p\symb}} = 1.
\end{equation}

By (B2) and since $\br_0(0) = 0$, for $\x \in \ui$, 
\begin{equation}\label{eq:decreasing-first-br}
  \ibr_0(\x) < \x,
\end{equation}
for all $\x \in \ui$. Let us now establish some useful equivalences:
\begin{align}
  \label{eq:ratio}
  & \frac{z_{p\symb} (\ibr_0 (\x))} {z_{p\symb} (\x)} \le \frac{\ell_{p\symb} (\ibr_0 (\x))}
  {\ell_{p\symb} (\x)} \\[3pt]
  \Longleftrightarrow \quad & \frac{z_{p\symb} (\x) + \int_{\x}^{\ibr_0(\x)} 
  z_{p\symb}'(s) \, ds} {z_{p\symb}(\x)} \le \frac{\ell_{p\symb}(\x) + 
  \int_{\x}^{\ibr_0(\x)} \ell_{p\symb}'(s) \, ds}{\ell_{p\symb}(\x)} \\
  \label{eq:ratio3}
  \Longleftrightarrow \quad & \ell_{p\symb}(\x) \int_{\ibr_0(\x)}^{\x} |z_{p\symb}'(s)| 
  \, ds \le z_{p\symb}(\x) \int_{\ibr_0(\x)}^{\x} |\ell_{p\symb}'(s)| \, ds.
\end{align}
So, (\ref{eq:ratio}) holds true by virtue of (\ref{eq:ratio3}), (\ref{eq:decreasing-first-br}) 
and (\ref{eq:mod-properties}). Therefore, using (\ref{eq:mannaggia}),
\begin{equation}
  \ibr_0' (\x) \, \frac{z_{p\symb} (\ibr_0 (\x))} {z_{p\symb} (\x)} \le \ibr_0' (\x) \, 
  \frac{\ell_{p\symb} (\ibr_0 (\x))} {\ell_{p\symb} (\x)} < 1,
\end{equation}
for all $\x \in (0,\delta)$. It remains to prove that the leftmost term above is less 
than 1 also for $\x \in [\delta,1]$. But this is an easy consequence of 
(\ref{eq:finefine}) and $\ibr_0' < 1$.
\qed

\begin{lemma} \label{lem:dist-int-2}
  For all $j$,
  \begin{displaymath}
    \sup_\ui \, \left| \frac{ \ds \frac{z'_{p\symb}} {z^2_{p\symb}} - \ibr'_j \, 
    \frac{ z'_{p\symb} \circ \ibr_j} {z_{p\symb} (z_{p\symb} \circ \ibr_j) }}
    {\ds 1 - | \ibr'_j | \, \frac{z_{p\symb} \circ \ibr_j} {z_{p\symb} }} \right| < \infty .
  \end{displaymath}
\end{lemma}

\proof Let us first consider the case $j=0$. Since $z_{p\symb} = \ell_{p\symb}$ in 
$(0, \eta/2)$, by Lemma \ref{lem:distortion-int} there exists $\eta_o \in (0,\eta/2)$ 
such that 
\begin{equation} \label{eq:boh}
\begin{split}
  & \sup_{(0,\eta_o)} \, \left| \frac{ \ds \frac{z'_{p\symb}} {z^2_{p\symb}} - \ibr'_0 \, 
  \frac{ z'_{p\symb} \circ \ibr_0} {z_{p\symb} (z_{p\symb} \circ \ibr_0) }}
  {\ds 1 - \ibr'_0 \, \frac{z_{p\symb} \circ \ibr_j} {z_{p\symb} }} \right| 
  = \sup_{(0,\eta_o)} \, \left| \frac{ \ds \frac{\ell'_{p\symb}} {\ell^2_{p\symb}} - \ibr'_0 \, 
  \frac{ \ell'_{p\symb} \circ \ibr_0} {\ell_{p\symb} (\ell_{p\symb} \circ \ibr_0) }}
  {\ds 1 - \ibr'_0 \, \frac{\ell_{p\symb} \circ \ibr_j} {\ell_{p\symb} }} \right| \\
  & \qquad \le \sup_{(0,\eta_o)} \frac{|\ell'_{p\symb}|} {\ell^2_{p\symb}} \, \cdot 
  \sup_{(0,\eta_o)} \, \frac{ 1-\ibr'_0 \, \ds \frac{ \ell_{p\symb} 
  (\ell'_{p\symb}\circ \ibr_0) } { \ell'_{p\symb} (\ell_{p\symb}\circ \ibr_0) } } 
  {1 - \ibr'_0 \, \ds \frac{ \ell_{p\symb} \circ \ibr_0} {\ell_{p\symb}} } < \infty 
\end{split}
\end{equation}
On the other hand, by Lemma \ref{lem:finite-modification}, using continuity, 
\begin{equation} \label{eq:boh2}
  \inf_{[\eta_o, 1]} \left(1 - \ibr'_0 \, \frac{z_{p\symb} \circ \ibr_0}{z_{p\symb}}
  \right) > 0.
\end{equation}
Since the numerator of the leftmost term in (\ref{eq:boh}) is bounded outside a 
neighborhood of the origin (recall that $z_{p\symb} \ge c > 0$), we have proved 
the claim of the Lemma when $j=0$. 

Let us now consider the case $j \ge 1$. Lemma \ref{lem:finite-modification} implies 
that 
\begin{equation}
  1 - | \ibr'_j | \, \frac{z_{p\symb}\circ \ibr_j}{z_{p\symb}} \ge 1 - \Lambda^{-1}>0.
\end{equation}
Also, $z'_{p\symb} \circ \ibr_j = 0$, because $\ibr_j \ge a_1 > \eta$, Thus,
\begin{equation}
  \sup_\ui \left| \frac{z'_{p\symb}} {z^2_{p\symb}} - \ibr'_j \, 
  \frac{ z'_{p\symb} \circ \ibr_j} {z_{p\symb} (z_{p\symb} \circ \ibr_j) } \right| = 
  \sup_\ui \frac{|z'_{p\symb}|}{z^2_{p\symb}} < \infty ,
\end{equation}
because the argument of the sup  is bounded in a right neighborhood of $0$
and identically null in a left neighborhood of 1.
The previous two estimates conclude the proof of Lemma \ref{lem:dist-int-2}.
\qed

\begin{lemma} \label{lem:dist-int-3}
  For all $j$,
  \begin{displaymath}
    \sup_\ui \, \frac{ \left| \ds \frac{\ibr_j''} {\ibr_j'} \right| } 
    { z_{p\symb} \left( \ds 1 - | \ibr'_j | \, \frac{z_{p\symb} \circ \ibr_j} 
    {z_{p\symb} } \right) } < \infty .
  \end{displaymath}
\end{lemma}

\proof We proceed as in the previous proof. For $j=0$, by Lemma 
\ref{lem:another-dist} there exists $\eta_o \in (0,\eta/2)$ such that
\begin{equation}
\begin{split}
  & \sup_{(0,\eta_o)}  \frac{ \ds \frac{ |\ibr_0'' | } {\ibr_0'} } 
  { z_{p\symb} \left( \ds 1 - \ibr'_0 \, \frac{z_{p\symb} \circ \ibr_0} {z_{p\symb} } 
  \right) } = \sup_{(0,\eta_o)} \frac{ \ds \frac{ |\ibr_0'' | } {\ibr_0'} } 
  { \ell_{p\symb} \left( \ds 1 - \ibr'_0 \, \frac{ \ell_{p\symb} \circ \ibr_0} 
  {\ell_{p\symb} } \right) } \\[4pt]
  & \qquad \le \frac{ \ds \sup_{(0,\eta_o)} \frac{ |\ibr_0'' | } {\ibr_0'}  }
  { \ds \inf_{(0,\eta_o)}  \ell_{p\symb} \left( \ds 1 - \ibr'_0 \, 
  \frac{ \ell_{p\symb} \circ \ibr_0} {\ell_{p\symb} } \right) } < \infty 
\end{split}
\end{equation}
(recall that is $\ibr_0''$ continuous on $[0,1]$). The bound on $[\eta_o, 1]$ 
follows from (\ref{eq:boh2}) and the definition of $z_{p\symb}$.

In the case $j \ge 1$, the assertion follows directly form the second statement 
of Lemma \ref{lem:finite-modification}, the definition of $z_{p\symb}$ and the fact 
that $|\ibr_j''|/|\ibr_j'|$ is bounded. 
\qed

As in the proof of Theorem \ref{thm-main-rp}, we need cones of local \ob s
and a suitable transfer operator to act on them.

\paragraph{Cones of local \ob s:} For $\delta, D, \eps > 0$, set
\begin{equation} \label{cones-ui}
  \mathcal C_{\delta,D,\eps} := \rset{g \in L^1(\zeta_{p\symb}) \cap C^1} 
  {g>0, \frac{|g'|}{g} \le D z_{p\symb} \mbox{ and } \frac{|g'(\x)|}{g(\x)} \le \eps 
  z_{p\symb}(\x), \forall \x \le \delta},
\end{equation}
where we have dropped from the indication that all \fn s are 
defined in $\ui$. These cones are analogous to the ones defined in 
(\ref{cones-rp}) for the half-line case, except that they are relative to the
distance 
\begin{equation}
  \mathrm{dist}(\x_1, \x_2) := \left| \int_{\x_1}^{\x_2} z_{p\symb}(\x) \, d\x \right|
\end{equation}
in $\ui$ (notice that the diameter of $\ui$ is infinity, w.r.t.\ dist).

\paragraph{Transfer operator:} We define $\trop_{p\symb}: L^1(\zeta_{p\symb}) 
\into L^1(\zeta_{p\symb})$ via the identity
\begin{equation} 
  \forall F \in L^{\infty}(\zeta_{p\symb}), g \in L^1(\zeta_{p\symb}), \qquad 
  \int_\ui (F \circ T) g \, d\zeta_{p\symb} = \int_\ui F \, \trop_{p\symb} g \, 
  d\zeta_{p\symb}.
\end{equation}
This operator describes the evolution of densities w.r.t.\ the \me\ 
$\zeta_{p\symb}$. Since $\frac{d\zeta_{p\symb}} {d\leb} = z_{p\symb}$, a 
standard computation gives
\begin{equation} \label{def-trop-psymb}
  \trop_{p\symb} \, g  = \frac{\trop (g z_{p\symb})}{z_{p\symb}} = \frac 1 {z_{p\symb}}
  \sum_{j=0}^{N-1} |\ibr'_j| \, (g z_{p\symb})\circ \ibr_j,
\end{equation}
where $\trop$ is the Perron-Frobenius operator for $T$ on $\ui$. 

Set
\begin{equation}
  \underline{D} := \max_{0 \le j \le N-1} \, \sup_\ui \left( \frac{ \ds 
  \frac{|\ibr''_j|}{|\ibr'_j|}} {z_{p\symb} \left( \ds 1 -  |\ibr'_j| \frac{z_{p\symb} \circ \ibr_j}
  {z_{p\symb}} \right) } + \frac{ \ds \left |\frac{z'_{p\symb}} {z_{p\symb}^2} - \ibr'_j \,
  \frac{z'_{p\symb} \circ \ibr_j} {z_{p\symb} (z_{p\symb} \circ \ibr_j)} \right |}
  {\left( \ds 1 -  |\ibr'_j| \, \frac{z_{p\symb} \circ \ibr_j} {z_{p\symb}} \right)}  \right).
\end{equation}
By Lemmas \ref{lem:dist-int-3} and \ref{lem:dist-int-2}, $\underline{D} < \infty$. 

\begin{lemma} \label{lem:cone-invariance-int}
  Let $D \ge \underline{D}$. For all $\eps >0$, there exists $\delta \in (0,1)$ 
  such that $\trop_{p\symb} \, \mathcal C_{\delta, D, \eps} \subseteq 
  \mathcal C_{\delta, D, \eps}$.
\end{lemma}

\proof A somewhat lengthy computation gives
\begin{equation} \label{eq:derivative-trans-interval}
  (\trop_{p\symb} \, g)' = \sum_{j=0}^{N-1} |\ibr'_j| \left[(g z_{p\symb}) \circ \ibr_j 
  \right] \mathcal D_j(g),
\end{equation}
where
\begin{equation} \label{djg}
  \mathcal D_j(g) := \frac1 {z_{p\symb}} \left( \frac{\ibr''_j}{\ibr'_j} + 
   \ibr'_j\frac{g' \circ \ibr_j }{g \circ \ibr_j} + 
  \ibr'_j\frac{z'_{p\symb}\circ \ibr_j}{z_{p\symb} \circ \ibr_j} - 
  \frac{z'_{p\symb}}{z_{p\symb}} \right),
\end{equation}
intended as a \fn\ $\ui \into \R$. For all $g>0$ with $\frac{|g'|}{g} \le D z_{p\symb}$, 
we have
\begin{equation} \label{eq:dist-int-cone}
  \max_j \|\mathcal D_j(g)\|_\infty \le \max_j \, \sup_\ui \left[ \frac1 {z_{p\symb}}
  \left( \left| \frac{\ibr''_j}{\ibr'_j} \right| + D (z_{p\symb} \circ \ibr_j) |\ibr'_j| +
  \left| \ibr'_j\frac{z'_{p\symb} \circ \ibr_j} {z_{p\symb} \circ \ibr_j}  
  - \frac{z'_{p\symb}}{z_{p\symb}} \right| \right) \right].
\end{equation}
Using that $D \ge \underline{D}$, one has that (\ref{eq:dist-int-cone}) is less than 
or equal to $D$. Therefore
\begin{equation} \label{eq:dist-int-cone2}
  \left| (\trop_{p\symb} \, g)' \right| \le \max_j \|\mathcal D_j(g) \|_{\infty} \,
  \trop (g z_{p\symb}) \le D \, \trop (g z_{p\symb}) = D z_{p\symb} \, 
  \trop_{p\symb} \, g,
\end{equation}
cf.\ (\ref{def-trop-psymb}).

For a given $\eps > 0$, we now describe how to choose $\delta$ so that 
$\mathcal C_{\delta, D, \eps}$ is invariant. There exists $\delta_1 \in 
(0, \eta/2)$, cf.\ (\ref{eq:modification-def}), such that the following 
functional (in)equalities are true, when restricted to $(0, \delta_1]$:
\begin{equation} \label{eq:delta-choice1}
\begin{split}
  & \left( \frac{\left| \ibr''_0 \right|}{z_{p\symb} \ibr'_0} + \left| \frac{z'_{p\symb}}
  {z_{p\symb}^2} - \ibr'_0 \, \frac{ (z'_{p\symb} \circ \ibr_0)} 
  {(z_{p\symb} \circ \ibr_0) z_{p\symb}} \right| \right) \left( 1 - \ibr'_0 \,\ds 
  \frac{z_{p\symb} \circ \ibr_0}{z_{p\symb}} \right)^{-1} \\[5pt]
  &\qquad = \frac{|\ibr''_0|}{ \ibr'_0 \, \ell_{p\symb} \left(1 - \ibr'_0 \,\ds 
  \frac{\ell_{p\symb} \circ \ibr_0}{\ell_{p\symb}} \right)} + \frac{|\ell'_{p\symb}|}
  {\ell_{p\symb}^2} \frac{ \left|1 - \ibr'_0 \,\ds \frac{\ell_{p\symb} 
  (\ell'_{p\symb} \circ \ibr_0)} {\ell'_{p\symb} (\ell_{p\symb} \circ \ibr_0)} \right|}
  {1 - \ibr'_0 \,\ds \frac{\ell_{p\symb} \circ \ibr_0}{\ell_{p\symb}}} \\[5pt]
  &\qquad \le \eps.
\end{split}
\end{equation}
In fact, calling $\x$ the argument of all the above \fn s, when $\x \to 0^+$
the second term of the second line of (\ref{eq:delta-choice1}) vanishes
by Lemma \ref{lem:distortion-int} and the fact that $\frac{|\ell'_{p\symb}(\x)|}
{\ell^2_{p\symb}(\x)} \to 0$. The same is true for the first term, as a 
consequence of Lemmas \ref{lem:second-derivative} and \ref{lem:another-dist}. 
Moreover, there exists $\delta_2 > 0$ such that, in $(0, \delta_2]$,
\begin{equation} \label{eq:delta-choice2}
  \frac{|\ibr''_j|}{z_{p\symb} |\ibr'_j|} + \frac{D (z_{p\symb} \circ \ibr_j) |\ibr'_j|}
  {z_{p\symb}} +  \frac{|z'_{p\symb}|}{z_{p\symb}^2} \le \eps,
\end{equation}
whenever $j \ge 1$. In fact, for the first two terms, both  $D |\ibr_j'| (z_{p\symb} 
\circ \ibr_j)$ and $|\ibr''_j / \ibr'_j|$ are bounded, for $j \ge 1$. As for the third
term, 
\begin{equation}
  \lim_{\x \to 0^+} \frac{ |z'_{p\symb}(\x)|}{z^2_{p\symb}(\x)} = \lim_{\x \to 0^+}
  \frac{|\ell'_{p\symb}(\x)|}{\ell^2_{p\symb}(\x)} = 0.
\end{equation}

Set $\delta := \min \{\delta_1, \delta_2\}$ and let $g \in 
\mathcal C_{\delta, D, \eps}$. We show that $\trop_{p\symb} \,g \in 
\mathcal C_{\delta, D, \eps}$. Evidently, $\trop_{p\symb} \,g > 0$ and 
$\trop_{p\symb}\, g \in C^1$. By (\ref{eq:dist-int-cone2}), 
$\frac{|\trop_{p\symb} \, g'|}{\trop_{p\symb} g} \le D z_{p\symb}$. For $\x \le \delta$, 
by (\ref{eq:decreasing-first-br}), $\phi_0(\x) \le \delta$, whence 
$\frac{|g'\circ \ibr_0(\x)|} {g \circ \ibr_0(\x)} \le \eps z_{p\symb} \circ \ibr_0(\x)$. 

Therefore, looking at (\ref{eq:derivative-trans-interval})-(\ref{djg}) and using that 
$z'_{p\symb} \circ \ibr_j \equiv 0$, for all $j \ge 1$, we can write that, on 
$(0,\delta]$, 
\begin{equation} \label{eq:derivative-small-int}
\begin{split}
  \left| (\trop_{p\symb} \, g)' \right| &\le \max_{j \ge 0} \left| \mathcal D_j(g) \right| 
  \trop (g z_{p\symb}) \\[2pt]
  &\le \max \left\{ \frac1 {z_{p\symb}} \left( \frac{|\ibr''_{0}|}{\ibr'_{0}} + \ibr'_0
  \, \eps (z_{p\symb} \circ \ibr_0) + \left| \frac{z'_{p\symb}}{z_{p\symb}} - 
  \frac{(z'_{p\symb}\circ \ibr_0) \ibr'_0} {z_{p\symb}\circ \ibr_0} \right| \right), 
  \right. \\[4pt]
  &\qquad \left. \max_{j \ge 1} \left[ \frac1 {z_{p\symb}} \left( \frac{|\ibr''_j|}{|\ibr'_j|} 
  + |\ibr'_j| D (z_{p\symb} \circ \ibr_j) + \frac{|z'_{p\symb}|}{z_{p\symb}} \right) 
  \right] \right\} \trop (g z_{p\symb}).
\end{split}
\end{equation}
Once again, the above inequality holds when all the \fn s are restricted to $(0,\delta]$. 
Using (\ref{eq:delta-choice1}) and (\ref{eq:delta-choice2}), one shows that, with the
same restriction on the arguments, both terms inside the above braces are 
$\le \eps$, concluding that, for $0 < \x \le \delta$,
\begin{equation}
  \left| (\trop_{p\symb}\, g)' (\x)\right|  \le \eps \trop(g z_{p\symb}) (\x) = \eps 
  z_{p\symb}(\x) (\trop_{p\symb} \,g)(\x) .
\end{equation}
This ends the proof of Lemma \ref{lem:cone-invariance-int}.
\qed

The analogue of Lemma \ref{lem:control-avg} is the following lemma, whose 
proof is included for the reader's convenience.

\begin{lemma} \label{lem:control-avg-int}
  Let $F \in \gou(\zeta_{p\symb})$, $\overline{\zeta_{p\symb}}(F) = 0$. For each 
  $\rho >0$, there  exists $\eps > 0$ such that for, every $g \in C^1(\ui)$ and 
  $\delta \ge 0$ with 
  \begin{itemize}
    \item[(i)] $g >0$;
    \item[(ii)] $\zeta_{p\symb}(g) = 1$;
    \item[(iii)] $\ds \frac{|g'(\x)|}{g(\x)} \le \eps z_{p\symb}(\x)$, for all $\x \le \delta$,
   \end{itemize} 
  one has
  \begin{displaymath}
    \left |\int_0^\delta F g \, d\zeta_{p\symb} \right| \le \rho.
  \end{displaymath}
\end{lemma}

\proof We assume $\|F\|_\infty > 0$, otherwise the statement is trivial. By the 
hypothesis $\overline{\zeta_{p\symb}}(F) = 0$, for any $\rho > 0$, there exists $r>0$
such that
\begin{equation} \label{eq:globavg-lem-int}
  \sup_{V \in \scv_{\zeta_{p\symb}}(r)} \left| \frac1r \int_V F \, d \zeta_{p\symb} \right| 
  \le \frac{\rho}2. 
\end{equation}
Set $\eps := \frac{\rho} {2r\|F\|_\infty}$. We partition $(0,\delta]$ into a countable
family of adjacent intervals $V_k \in \scv_{\zeta_{p\symb}}(r)$. For all $g$ and 
$\delta$ verifying \emph{(i)-(iii)}, 
\begin{equation} \label{eq:split-integral-int}
  \int_0^\delta F g \, d \zeta_{p\symb} = \sum_{k=0}^\infty \int_{V_k} F g \, 
  d \zeta_{p\symb}.
\end{equation}
For $k \in \N$, denote by $\x_{\mathrm{min},k}$ a minimum point of $g$ in $V_k$. 
By (\ref{eq:globavg-lem-int}), \emph{(i)}, \emph{(iii)} and the definition of $\eps$,
\begin{equation} \label{eq:small-pieces-int}
\begin{split}
  \left| \int_{V_k} F g \, d \zeta_{p\symb} \right| &= 
  \left| \int_{V_k} F(\x) \left( g(\x_{\mathrm{min},k}) + 
  \int_{\x_{\mathrm{min},k}}^\x g'(t) \, dt \right) d \zeta_{p\symb}(\x) \right| \\
  &\le \frac \rho 2 r g(\x_{\mathrm{min},k}) + r \|F\|_\infty \, \eps 
  \left(\int_{V_k} g(t) \, z_{p\symb}(t) \, dt \right) \\
  &\le \left(\frac \rho 2 + \frac \rho 2 \right) \left(\int_{V_k} g(t) \, 
  d \zeta_{p\symb}(t) \right).
\end{split}
\end{equation}
and we conclude the proof of the lemma by (\ref{eq:split-integral-int}), 
(\ref{eq:small-pieces-int}) and \emph{(ii)}.
\qed

As a penultimate step towards the proof of Theorem \ref{thm-main-ui}, we show
that $T$ is fully global-local mixing w.r.t.\ $\scv$ and the \me\ $\zeta_{p\symb}$. 
To this end, in complete analogy with the case of maps of $\rp$ --- see the 
discussion before (\ref{goal2}) --- it suffices to prove that, for all $F \in 
\gou(\zeta_{p\symb})$ with $\overline{\zeta_{p\symb}}(F) = 0$ and $\rho>0$, 
there exists $g_\rho \in L^1(\zeta_{p\symb})$ with $\zeta_{p\symb}(g_{\rho}) = 1$ 
such that
\begin{equation} \label{eq:goal2-in}
  \limsup_{n \to \infty} \, | \zeta_{p\symb} ((F\circ T^n)g_{\rho}) | \le 2\rho.
\end{equation}
To this end, given $\rho>0$, we consider $\eps$ as given by Lemma 
\ref{lem:control-avg-int} and $\delta$ as given by Lemma 
\ref{lem:cone-invariance-int} for such value of $\eps$. Let $D \ge \underline{D}$ 
and $g_\rho$ be any element of $\mathcal C_{\delta,D,\eps}$ with 
$\zeta_{p\symb}(g_\rho) =1$. For all $n \in \N$,
\begin{equation}
\begin{split}
  | \zeta_{p\symb}((F \circ T^n) g_\rho) | &\le | \zeta_{p\symb}( (1_{[\delta,1]} F) 
  \circ T^n) g_\rho) | + | \zeta_{p\symb}( (1_{(0, \delta)} F) \circ T^n) g_\rho) | \\
  &= | \zeta_{p\symb}( (1_{[\delta,1]} F) \circ T^n) g_\rho) | + 
  \left| \int_0^\delta F (\trop_{p\symb}^n \, g_\rho) \, d\zeta_{p\symb} \right| .
\end{split}
\end{equation}
For $n \to \infty$, the first term vanishes by local-local mixing, for $T$ is exact
(w.r.t\ $z_{p\symb}$ as well as $\leb$). The second term is bounded above by 
$\rho$ by Lemma \ref{lem:control-avg-int}, since $\trop_{p\symb}^n \, g_\rho \in 
\mathcal C_{\delta,D,\eps}$  (Lemma \ref{lem:cone-invariance-int}) and 
$\zeta_{p\symb} (\trop_{p\symb}^n \, g_\rho) = \| \trop_{p\symb}^n \, 
g_\rho \|_{L^1(\zeta_{p\symb})} =  \| g_\rho \|_{L^1(\zeta_{p\symb})} = 1$.

Finally, full global-local mixing w.r.t.\ $\zeta_{p\symb}$ and w.r.t.\ 
$\lambda_{p\symb}$ are equivalent properties by Proposition \ref{prop2} 
and Lemma \ref{lem-primo-ultimo}. This proves the statement of Theorem 
\ref{thm-main-ui} for the case $p>1$. As for the case $p=1$, one sees by 
Proposition \ref{prop:restriction} and Remark \ref{rm:simpl-cond} that 
$\overline{\zeta_{1+}}$ is a restriction of $\overline{\lambda_1}$. But
$\overline{\lambda_{1+}} = \overline{\zeta_{1+}}$ (Proposition \ref{prop2}),
so $\overline{\lambda_{1+}}$ is a restriction of $\overline{\lambda_1}$.%
\footnote{A strict restriction, in fact, as one can see that $F_1(\x) = 
\cos( \log(\x^{-1}))$ belongs to  $\gou(\lambda_1)$ but not to 
$\gou(\lambda_{1+})$; cf.\ proof of Corollary \ref{cor-main-ui} in Section 
\ref{subs-maps-ui}.}
Then one applies Lemma \ref{lem-primo-ultimo} 
to conclude that $(\ui, \scb, \lambda_1, T)$ is global-local \m\ 
w.r.t.\ $\scv$, $\gou(\lambda_{1+})$ and $L^1(\lambda_1)$.
\qed

\footnotesize

\end{document}